\documentclass[11pt]{amsart}

\usepackage{amscd,array,graphicx,color,epsfig}

\newtheorem{theorem}{Theorem}[section]
\newtheorem{proposition}[theorem]{Proposition}
\newtheorem{lemma}[theorem]{Lemma}
\newtheorem{corollary}[theorem]{Corollary}
\newtheorem{parameterizationtheorem}[theorem]{Parameterization Theorem}
\newtheorem*{ParameterizationTheorem}{Parameterization Theorem~\ref{thm:parameterization}}
\newtheorem{uniquecablingsequencetheorem}[theorem]{Unique Cabling Sequence Theorem}
\newtheorem*{UniqueCablingSequenceTheorem}{Unique Cabling Sequence Theorem~\ref{thm:cable_sequence}}

\newtheorem*{MinimumBridgeNumberTheorem}{Minimum Bridge Number Theorem}

\theoremstyle{definition}
\newtheorem{definition}[theorem]{Definition}
\newtheorem{notation}[theorem]{Notation}

\newtheorem{remark}[theorem]{Remark}

\numberwithin{equation}{section}

\newcommand{\longpage}{\enlargethispage{\baselineskip}}
\newcommand{\shortpage}{\enlargethispage{-\baselineskip}}

\newcommand{\Diff}{\operatorname{Diff}}
\newcommand{\depth}{\operatorname{depth}}

\newcommand{\lk}{\operatorname{lk}}

\newcommand{\Out}{\operatorname{Out}}
\newcommand{\PSL}{\operatorname{PSL}}
\newcommand{\SL}{\operatorname{SL}}

\newcommand{\Q}{\operatorname{{\mathbb Q}}}

\newcommand{\Z}{\operatorname{{\mathbb Z}}}

\newcommand{\D}{\operatorname{{\mathcal D}}}
\newcommand{\G}{\operatorname{{\mathcal G}}}
\newcommand{\K}{\operatorname{{\mathcal K}}}
\newcommand{\T}{\operatorname{{\mathcal T}}}

\begin{document}

\title[The tree of knot tunnels]
{The tree of knot tunnels}

\author{Sangbum Cho}
\address{Department of Mathematics\\
University of Oklahoma\\
Norman, Oklahoma 73019\\
USA} 
\email{scho@ou.edu}

\author{Darryl McCullough}
\address{Department of Mathematics\\
University of Oklahoma\\
Norman, Oklahoma 73019\\
USA} 
\email{dmccullough@math.ou.edu}
\urladdr{www.math.ou.edu/$_{\widetilde{\phantom{n}}}$dmccullough/}

\subjclass{Primary 57M25}

\date{\today}

\keywords{knot, link, tunnel, (1,1), disk complex, two-bridge}

\begin{abstract}
We present a new theory which describes the collection of all tunnels of
tunnel number $1$ knots in $S^3$ (up to orientation-preserving equivalence
in the sense of Heegaard splittings) using the disk complex of the
genus-$2$ handlebody and associated structures. It shows that each knot
tunnel is obtained from the tunnel of the trivial knot by a uniquely
determined sequence of simple cabling constructions. A cabling construction
is determined by a single rational parameter, so there is a corresponding
numerical parameterization of all tunnels by sequences of such parameters
and some additional data. Up to superficial differences in definition, the
final parameter of this sequence is the Scharlemann-Thompson invariant of
the tunnel, and the other parameters are the Scharlemann-Thompson
invariants of the intermediate tunnels produced by the constructions.  We
calculate the parameter sequences for tunnels of $2$-bridge knots.  The
theory extends easily to links, and to allow equivalence of tunnels by
homeomorphisms that may be orientation-reversing.
\end{abstract}

\maketitle

\section*{Introduction}
\label{sec:intro}

In this work we present a new descriptive theory for the tunnels of tunnel
number~$1$ knots in $S^3$. At its heart is a bijective correspondence
between the set of equivalence classes of all tunnels of all tunnel number
1 knots and a subset of the vertices of a certain tree $\T$. In fact, $\T$
is bipartite, and the tunnel vertex subset is exactly one of its two
classes of vertices. The construction of $\T$ uses the disk complex of the
genus~$2$ handlebody, and $\T$ is a quotient of a spine of the subcomplex
of nonseparating disks. The tree and its associated objects have a rich
combinatorial structure. The work in this paper is a first step toward
understanding how that structure is manifested in the topology of tunnel
number $1$ knots.

The theory has many consequences. It shows that every tunnel can be
obtained by starting from the unique tunnel of the trivial knot and
performing a uniquely determined sequence of simple constructions. Each
construction is determined by a rational parameter, which is essentially a
Scharlemann-Thompson invariant \cite{Scharlemann-Thompson}. This gives a
natural numerical parameterization of all tunnels. We have computed this
sequence of rational invariants for all tunnels of $2$-bridge knots, and
for the ``short'' tunnels of torus knots.

The theory adapts easily to include tunnels of tunnel number $1$ links, and
to allow orientation-reversing equivalence.

The next two sections provide overviews before beginning the actual
development of the theory. Section~\ref{sec:context} emphasizes the context
of the work, while section~\ref{sec:overview} summarizes the main ideas and
results of the paper and the contents of the individual sections.
\longpage

\section{Context of the work}
\label{sec:context}

There are several equivalent definitions of a tunnel for a knot $K$ in
$S^3$. One is that a tunnel is a $1$-handle attached to a regular
neighborhood of $K$ in $S^3$ to produce a genus-$2$ handlebody that is
unknotted, that is, one of the handlebodies of a genus-$2$ Heegaard
splitting of $S^3$. Two such configurations are equivalent tunnels when the
handlebodies are isotopic taking the copy of $K$ in one handlebody to the
corresponding copy in the other. Alternatively, one may think of a tunnel
as an arc $\alpha$ meeting $K$ only in its endpoints, such that a regular
neighborhood of the ``$\theta$-curve'' $K\cup\alpha$ is unknotted. The
equivalence on such arcs must then allow not only isotopy but also
``sliding,'' where the two endpoints can meet and pass through each other.

There is a stronger notion of equivalence, in which the isotopies and
sliding must preserve the knot at all times. All of our work uses the
weaker notion.

If one removes a small tubular neighborhood of $K$ from the genus-$2$
handlebody produced by a tunnel, the remaining compression body and its
complementary genus-$2$ handlebody form a genus-$2$ Heegaard splitting of
the knot space, in the sense of Heegaard splittings of $3$-manifolds with
boundary. For closed $3$-manifolds, genus-$0$ Heegaard splittings are
trivial, and genus-1 splittings are very restrictive, forming only lens
spaces (including $S^3$ and $S^2\times S^1$), while genus-$2$ splittings
are already a very complicated class. In this context, the Heegaard
splittings coming from tunnel number 1 knots might be considered to be a
special class of ``genus-$1\frac12$'' splittings, an intermediate case
where one might hope to find structure restricted enough to be tractable,
but rich enough to be of mathematical interest.

The historical development of the subject is consistent with this hope.  An
impressive amount of geometric theory of tunnel number $1$ knots has been
developed by a number of researchers. Recently, a general picture has begun
to emerge, through work of M. Scharlemann and A. Thompson
\cite{Scharlemann-Thompson} which defines a rational invariant that detects
a kind of cabled structure of $K$ near the tunnel. One of the applications
of our work is a complete clarification of how their invariant works and
what information it is detecting. As mentioned above, it extends to a
sequence of rational invariants which describe a unique sequence of simple
constructions that produce the tunnel. These rational invariants, plus a
bit more information, give a natural numerical parameterization of all knot
tunnels.

Other recent work in the subject has begun to utilize connections between
tunnel number $1$ knots and the curve complex of the genus-$2$ surface.
The curve complex provides an important measure of complexity of Heegaard
splittings, called the \textit{(Hempel) distance.} Applying it to the
splittings that correspond to knot tunnels, J. Johnson and A. Thompson
\cite{JohnsonBridgeNumber,Johnson-Thompson} and Y. Minsky, Y. Moriah, and
S. Schleimer \cite{MMS} have obtained results on bridge number and other
aspects of tunnels. In~\cite{CM}, we use our theory to define a new
distance-type invariant for knot tunnels, the ``depth'', that is finer than
distance (depth can be large even when distance is small). As one
application of the depth invariant, we can substantially improve known
estimates of the growth rate of bridge number as distance increases.

There are additional reasons to believe that the class of tunnel number $1$
knots is a nexus of interesting mathematical objects. The fundamental group
of the genus-$2$ handlebody is the free group on two generators $F_2$, and
because of this, the mapping class group of the genus-$2$ handlebody is
related to the automorphism group $\Out(F_2)$ and thereby to the linear
groups $\SL_2$, which are very special and whose theory differs in many
respects from the $\SL_n$ with~$n\geq 3$. It will be evident that the tree
$\T$ that is the central object in our theory carries much of the structure
of the well-known tree associated to~$\PSL_2(\Z)$.

\section{Summary of the results}
\label{sec:overview}

In this section, we will give an overview of the theory and applications
developed in this paper, and our paper~\cite{CM}.

Let $H$ be an unknotted genus-$2$ handlebody in $S^3$. By $\D(H)$ we denote
the $2$-dimensional complex whose vertices are the isotopy classes of
\textit{nonseparating} essential properly-imbedded disks in $H$. A set of
vertices spans a simplex of $\D(H)$ when they have pairwise disjoint
representatives. In the first barycentric subdivision $\D'(H)$ of $\D(H)$,
the span of the vertices that are not vertices of $\D(H)$ is a tree
$\widetilde{\T}$. Each vertex of $\widetilde{\T}$ is either a triple of
(isotopy classes of) disks in $H$, or a pair. Figure~\ref{fig:subdivision}
below shows a small portion of $\D(H)$ and~$\widetilde{\T}$.

Each tunnel of a tunnel number 1 knot determines a collection of disks in
$H$ as follows. The tunnel is a $1$-handle attached to a regular
neighborhood of the knot to form an unknotted genus-$2$ handlebody. An
isotopy carrying this handlebody to $H$ carries a cocore $2$-disk of that
$1$-handle to a nonseparating disk in $H$. The indeterminacy of this
process is the group of isotopy classes of orientation-preserving
homeomorphisms of $S^3$ that preserve $H$. This is the \textit{Goeritz
group $\G$,} which we discuss in section~\ref{sec:Goeritz} below. In
particular, we will recall work of M. Scharlemann \cite{ScharlemannTree}
and E. Akbas \cite{Akbas} that shows that $\G$ is finitely presented, and
even provides a simple presentation of it. 

Since $\G$ is exactly the indeterminacy of the cocore disk, moved to $H$,
it is essential to understand the action of $\G$ on $\D(H)$. This action is
closely related to a central concept of our viewpoint, called
primitivity. A disk in $H$ is \textit{primitive} if there exists a properly
imbedded disk $\tau'$ in the complementary handlebody $\overline{S^3-H}$
such that the circles $\partial\tau$ and $\partial\tau'$ in $\partial H$
intersect transversely in a single point. Viewed as tunnels, primitive
disks are exactly the tunnels of the trivial knot, and are all equivalent.
The portions of $\D(H)$ and $\widetilde{\T}$ corresponding to primitive
disks form the ``primitive region,'' in particular the pairs and triples of
primitive disks span a $\G$-invariant subtree $\widetilde{\T}_0$ of
$\widetilde{\T}$ called the \textit{primitive subtree.} The primitive
subtree $\widetilde{\T}_0$ is isomorphic to the tree used by Scharlemann
and Akbas to understand the Goeritz group. As explained in
section~\ref{sec:Goeritz} below, the $2$-complex used by Scharlemann imbeds
in a very natural way into $\D(H)$, in such a way that it deformation
retracts to $\widetilde{\T}_0$ (see figure~\ref{fig:primitive_tree}).
Indeed, it is fair to say that $\widetilde{\T}_0$ \textit{is} the
Scharlemann-Akbas tree. This viewpoint leads to a new proof of the main
results of \cite{ScharlemannTree} and \cite{Akbas} that avoids many of the
difficult geometric arguments of those important papers. This recasting of
their work is carried out in~\cite{Cho}.

Using this viewpoint, it is not difficult to analyze the action of $\G$ on
$\D(H)$ and $\widetilde{\T}$. The action of $\G$ on the primitive
structures is as transitive as possible, while the action of $\G$ on the
nonprimitive structures has stabilizers that are as small as possible---
usually only the order-$2$ subgroup generated by the ``hyperelliptic''
element which acts trivially on all of $\D(H)$. This gives a picture of
$\D(H)/\G$ as a tiny primitive region, together with additional portions
which look exactly as they did in~$\D(H)$. Figure~\ref{fig:Delta} below
shows the structure of $\D(H)/\G$ and~$\widetilde{\T}/\G$ near the
primitive region.  The quotients $\D(H)/\G$ and $\widetilde{\T}/\G$ are
described precisely in section~\ref{sec:quotients}, after preliminary work
which provide a useful framework for that discussion and is used throughout
the rest of the paper as well.

The tree $\T$ discussed in the introduction is defined to be
$\widetilde{\T}/\G$. The vertices of $\D(H)/\G$ are not in $\T$, but their
links in the barycentric subdivision $\D'(H)/\G$ are subcomplexes of
$\T$. These links are infinite trees. In $\T$, there is a vertex $\theta_0$
which is the unique $\G$-orbit of a triple of primitive disks in $H$ (or
dually it represents a planar $\theta$-curve in $S^3$). For each vertex
$\tau$ of $\D(H)/\G$, i.~e.~each tunnel, there is a unique shortest path in
$\T$ from $\theta_0$ to \textit{the vertex in the link of $\tau$ that is
closest to $\theta_0$.}  This closest vertex is a triple, called the
\textit{principal vertex} of $\tau$, and the path is the \textit{principal
path} of~$\tau$. Figure~\ref{fig:principal_vertex} below illustrates these
structures.

Those familiar with \cite{Scharlemann-Thompson} will not be surprised to
learn that the disks in the principal vertex, other than $\tau$ itself, are
the disks $\mu^+$ and $\mu^-$ used in the definition of the
Scharlemann-Thompson invariant. That is, the principal vertex of $\tau$ is
$\{\mu^+,\mu^-,\tau\}$.

The principal path and principal vertex of $\tau$, along with the
surrounding structure of $\D(H)/\G$, encode a great deal of geometric
information about the tunnel $\tau$ and the knot $K_\tau$ of which it is a
tunnel. In section~\ref{sec:cabling_construction} we examine how the
sequence of vertices in the principal path corresponds to a sequence of
constructions of a very simple type, which will look familiar to
experts. In a sentence, the construction is ``Think of the union of $K$ and
the tunnel arc as a $\theta$-curve, and cable the ends of the tunnel arc
and one of the arcs of $K$ in a neighborhood of the other arc of $K$.''
Each of these ``cabling constructions'' (often just called ``cablings'') is
determined by a rational ``slope'' parameter. For the first cabling of the
sequence, the indeterminacy coming from the Goeritz group necessitates some
special treatment; in particular, its parameter is a ``simple'' slope
taking values in $\Q/\Z$.  Figures~\ref{fig:schematic}
and~\ref{fig:cabling} below should give a fairly good idea of how cabling
constructions change a pair consisting of a knot and tunnel to a more
complicated pair.

More precisely, the vertices in the principal path of $\tau$ are a sequence
of alternating triples and pairs, which we write as $\theta_0$, $\mu_0$,
$\mu_0\cup \{\tau_0\}$, $\mu_1,\ldots\,$, $\mu_n$, $\mu_n\cup \{\tau_n\}$,
where $\tau_n=\tau$.  Denoting by $(\mu;\sigma)$ the edge of $\T$ that goes
from a pair $\mu$ to an adjacent triple $\mu\cup{\sigma}$, the exact
information needed to describe a cabling construction is a succession from
$(\mu_{i-1};\tau_{i-1})$ to $(\mu_i;\tau_i)$ in the principal path. The
first edge $(\mu_0;\tau_0)$ is special, and by itself determines a
``simple'' cabling. The uniqueness of the principal path gives one of our
main results:
\begin{UniqueCablingSequenceTheorem}
Let $\tau$ be a tunnel of a nontrivial knot.  Let $\theta_0$, $\mu_0$,
$\mu_0\cup \{\tau_0\}$, $\mu_1,\ldots\,$, $\mu_n$, $\mu_n\cup \{\tau_n\}$
with $\tau_n=\tau$ be the principal path of $\tau$.  Then the sequence of
$n+1$ cablings consisting of the simple cabling determined by
$(\mu_0;\tau_0)$ and the cablings determined by the successions from
$(\mu_{i-1};\tau_{i-1})$ to $(\mu_i;\tau_i)$ is the \textit{unique}
sequence of cablings beginning with the tunnel of the trivial knot and
ending with $\tau$.
\end{UniqueCablingSequenceTheorem}

As detailed in sections~\ref{sec:slope_disks} and~\ref{sec:general_slopes},
an edge $(\mu;\sigma)$ defines a coordinate system in which each essential
disk in $H$ disjoint from the disks of $\mu$ and not parallel to either of
them or to $\sigma$ is assigned a rational slope. In particular, the
$(\mu_i;\sigma_i)$-slopes of the $\tau_i$, where $\sigma_i$ is the unique
disk in $\mu_{i-1}-\mu_i$, together with a special $\Q/\Z$-valued slope
associated to the initial $(\mu_0;\tau_0)$, determine the exact cabling
constructions in the sequence. This gives us a version of
theorem~\ref{thm:cable_sequence} that describes the unique cabling
sequences as a parameterization of all tunnels:
\begin{ParameterizationTheorem}
Let $\tau$ be a knot tunnel with principal path $\theta_0$, $\mu_0$,
$\mu_0\cup \{\tau_0\}$, $\mu_1,\ldots\,$, $\mu_n$, $\mu_n\cup \{\tau_n\}$.
Fix a lift of the principal path to $\D(H)$, so that each $\mu_i$
corresponds to an actual pair of disks in~$H$.
\begin{enumerate}
\item If $\tau$ is primitive, put $m_0=[0]\in\Q/\Z$. Otherwise,
let $m_0=[p_0/q_0]\in\Q/\Z$ be the simple slope of $\tau_0$.
\item If $n\geq 1$, then for $1\leq i\leq n$ let $\sigma_i$ be the unique disk
in $\mu_{i-1}-\mu_i$ and let $m_i=q_i/p_i\in\Q$ be the
$(\mu_i;\sigma_i)$-slope of $\tau_i$.
\item If $n\geq 2$, then for $2\leq i\leq n$ define $s_i=0$ or $s_i=1$
according to whether or not the unique disk of $\mu_i\cap\mu_{i-1}$ equals the
unique disk of $\mu_{i-1}\cap\mu_{i-2}$.
\end{enumerate}
Then, sending $\tau$ to the pair $((m_0,\ldots,m_n),(s_2,\ldots,s_n))$ is a
bijection from the set of all tunnels of all tunnel number~$1$ knots to the
set of all elements $(([p_0/q_0],q_1/p_1,\ldots,q_n/p_n),(s_2,\ldots,s_n))$
in
\[\big(\Q/\Z \big)\,\cup\, 
\big(\Q/\Z\,\times\, \Q\big) \,\cup\, \big(\cup_{n\geq 2}
\;\Q/\Z\,\times\, \Q^n \,\times\,\, C_2^{n-1}\big)\]
with all $q_i$ odd.
\end{ParameterizationTheorem}
\noindent This is actually proven earlier than
theorem~\ref{thm:cable_sequence}, since it does not require any
interpretation of the principal path in terms of cabling constructions.

Some examples should be mentioned. A tunnel produced from the tunnel of the
trivial knot by a single cabling construction is called a \textit{simple}
tunnel. These are exactly the ``upper and lower'' tunnels of $2$-bridge
knots. According to theorem~\ref{thm:parameterization}, these are
determined by a single $\Q/\Z$-valued parameter $m_0$, and this is of
course a version of the standard rational parameter associated to the
$2$-bridge knot. Simple tunnels are examined in
section~\ref{sec:simple_tunnels}.

A non-simple tunnel produced by a cabling sequence in which one of the
original arcs of the trivial knot is retained is called a ``semisimple''
tunnel. In terms of the parameterization, the simple and semisimple tunnels
are exactly the tunnels with all $s_i=0$. Geometrically, the simple and
semisimple tunels are the ``(1,1)'' tunnels of $(1,1)$-knots (i.~e.~knots
which can be put in $1$-bridge position with respect to the levels of a
product neighborhood of an unknotted torus in $S^3$). As arcs, they are
sometimes called ``eyeglass'' tunnels, meaning that the tunnel arc can be
slid to an unknotted circle. The ``semisimple region'' in $\D(H)/\G$
appears to be where the more complicated phenomena involving tunnels
occur. Indeed, we do not know an example of a knot with inequivalent
tunnels which is not a $(1,1)$-knot.

The tunnels of $2$-bridge knots are examined in section~\ref{sec:2bridge}.
It is known from work of several mathematicians \cite{Kobayashi1,
Kobayashi2, Morimoto-Sakuma, Uchida} that a $2$-bridge knot has at most
four equivalence classes of tunnels (not six, for us, since we are
considering tunnels only up to equivalence, rather than up to isotopy). Two
of these are the upper and lower simple tunnels. The others are semisimple
tunnels, and in section~\ref{sec:2bridge} we determine their exact cabling
sequences. Indeed, we have made software available~\cite{slopes} that
computes them quite effectively.

In \cite{CM}, we examine another kind of knot whose tunnels have been
classified, the torus knots. For most torus knots, there are two semisimple
tunnels, plus a third ``short'' tunnel represented by an arc cutting across
the complementary annulus in a torus containing the knot. From
\cite{Scharlemann-Thompson}, it is known that the slope invariants of the
short tunnel are integral, and we determine their exact values. Not
surprisingly, for a $(p,q)$ torus knot the calculation uses a continued
fraction expansion of~$p/q$. Software for this is also
available~\cite{slopes}.

Section~\ref{sec:STinvariant} explains how the Scharlemann-Thompson
invariant is really the slope parameter for the final cabling
construction--- the $m_n$ in theorem~\ref{thm:parameterization}. This
rational number is called the \textit{principal slope} of $\tau$.  Because
of differing definitions, the invariants generally have different values,
but they capture exactly the same geometric information. Intuitively, the
Scharlemann-Thompson invariant is ``the slope of the disk that the tunnel
disk replaced as seen from the tunnel disk'' while the principal slope is
``the slope of the tunnel disk as seen from the disk it replaced,'' so it
is not surprising that they are related by a continued fraction algorithm,
which is described in proposition~\ref{prop:change_coordinates} and has
been computationally implemented~\cite{slopes}. In particular, when one of
them is an integer, the other is also an integer, in fact, the negative of
the first one.

Our entire theory extends quite easily to links, as we discuss in
section~\ref{sec:links}. The cabling constructions are expanded to allow
cablings that produce links, which are terminal in the sense that they do
not allow further cabling constructions to be performed. The
Parameterization Theorem~\ref{thm:parameterization} holds as stated, except
allowing the final $q_n$ to be even. One application of the link version of
the theory is a very quick proof of the fact that the only tunnels of a
$2$-bridge link are its upper and lower tunnels~\cite{Adams-Reid,Kuhn}. We
also show that a tunnel number~$1$ link with an unknotted component must
have torus bridge number~$2$ (theorem~\ref{thm:semisimple_links}).

The minimal length of a simplicial path in the $1$-skeleton of $\D(H)/\G$
from a tunnel $\tau$ to the orbit of the primitive disks is called the
\textit{depth} of~$\tau$. The simple and semisimple tunnels are exactly the
tunnels of depth~$1$. The depth invariant is the subject of our
paper~\cite{CM}. In particular, we show there that as a rather immediate
consequence of work of Goda, Scharlemann, and Thompson \cite{G-S-T} and
Scharlemann and Thompson \cite{Scharlemann-Thompson}, the bridge number
grows exponentially with the depth. More precisely, we have the following
from~\cite{CM}:
\begin{MinimumBridgeNumberTheorem} 
For $d\geq 1$, the minimum bridge number of a knot having a tunnel of
depth~$d$ is $a_d$, where $a_1=2$, $a_2=4$, and $a_d=2a_{d-1}+a_{d-2}$ for
$d\geq 3$.
\end{MinimumBridgeNumberTheorem}
As a matrix, the recursion in the Minimum Bridge Number Theorem is
\[\begin{pmatrix} a_{d+1}\\
a_d\end{pmatrix} =
\begin{pmatrix} 2 & 1 \\
1 & 0
\end{pmatrix}
\begin{pmatrix} a_d\\
a_{d-1}
\end{pmatrix}\ .
\]
The eigenvalues of this matrix are $1\pm \sqrt{2}$,
showing that the asymptotic growth of the bridge numbers of any sequence of
tunnel number~$1$ knots as a function of depth is at least a constant
multiple of $(1+\sqrt{2})^d$. This improves Lemma~2 of
\cite{JohnsonBridgeNumber}, which is that bridge number grows linearly with
distance. It also improves Proposition~1.11 of \cite{G-S-T}, which shows
that bridge number grows asymptotically at least as fast as~$2^d$.

Of course, the Minimum Bridge Number Theorem also shows that the bound for
growth rate of $(1+\sqrt{2})^d$ is best possible, indeed its proof tells
how to construct a tunnel of depth~$d$ having bridge number~$a_d$.
In~\cite{CM}, it is also shown that minimum growth rate is achieved by a
sequence of torus knot tunnels. In fact, the bridge numbers of that
sequence are given by the recursion in the Minimum Bridge Number Theorem,
except that one starts with $a_1=2$ and $a_2=5$. The terms of this
sequence are the minimal bridge numbers of any torus knot having a tunnel
of the corresponding depth.

The depth invariant has a natural geometric interpretation in terms of a
construction that first appeared in the paper of H. Goda, M. Scharlemann,
and A. Thompson \cite{G-S-T}. That construction, which we call a giant
step, takes a tunnel and produces a new tunnel (usually of a different
knot). They proved that every tunnel could be produced starting from the
tunnel of the trivial knot and applying giant steps. It turns out that
$\depth(\tau)$ is the minimum length of a sequence of giant steps that
produces $\tau$. Unlike the construction of a knot tunnel using cabling
operations, the choice of giant steps is usually not unique, even when one
restricts to minimal sequences. Using the combinatorial structure of
$\D(H)/\Gamma$, we give in~\cite{CM} an algorithm to calculate the number
of distinct minimal sequences of giant steps that produce a given
tunnel. In particular, this provides arbitrarily complicated examples of
tunnels for which the minimal giant step construction \textit{is} unique,
while showing that such tunnels are sparse among the set of all tunnels.

Our entire theory also adapts easily to allow tunnel equivalences which may
be orientation-reversing homeomorphisms of $S^3$. In the Parameterization
Theorem~\ref{thm:parameterization}, the cabling sequence of the mirror
image of a tunnel $\tau$ has the same parameters as $\tau$ except that the
slopes $(m_0,\ldots,m_n)$ become $(-m_0,\ldots,-m_n)$. This shows
(theorem~\ref{thm:HopfLink}) that apart from the tunnels of the trivial
knot and trivial link, the only tunnel that is equivalent to itself under
an orientation-reversing homeomorphism of $S^3$ is the tunnel of the Hopf
link. In fact, the Hopf link and its tunnel have some unusual symmetries,
which are analyzed in section~\ref{sec:Hopf_link}.

\section{The disk complex of an irreducible $3$-manifold}

Let $M$ be a compact, irreducible $3$-manifold. The \textit{disk complex}
$\K(M)$ is the simplicial complex whose vertices are the isotopy classes of
essential properly imbedded disks in $M$, such that a collection of $k+1$
vertices spans a $k$-simplex if and only if they admit a set of
pairwise-disjoint representatives.

This is a good point at which to mention that to avoid endless repetition
of ``isotopy class'', we often speak of ``disks'' and other objects when we
really mean their isotopy classes, with the implicit understanding that we
always choose representatives of the isotopy classes that are the simplest
possible with respect to whatever we are doing (i.~e.~transversely
intersecting in the minimum possible number of components, and so on).  A
``unique'' disk means a unique isotopy class of disks. We often omit
other formalisms that should be obvious from context, so imbedded
submanifolds are assumed to be essentially and properly imbedded, unless
otherwise stated, and isotopies are assumed to preserve relevant
structure. To initiate this massive abuse of language, we say that a
collection of $k+1$ disks in $M$ spans a $k$-simplex of $\K(M)$ if and only
if they are pairwise disjoint.

The following was proven in~\cite[Theorem~5.3]{JDiff}:
\begin{theorem}
If $\partial M$ is compressible, then $\K(M)$ is contractible.
\end{theorem}

Denote by $\D(M)$ the subcomplex of $\K(M)$ spanned by the
\textit{nonseparating} disks. From~\cite[Theorem~5.4]{JDiff}, we have
\begin{theorem} If $M$ has a nonseparating compressing disk, then
$\D(M)$ is contractible.
\label{thm:contractible}
\end{theorem}
The basic idea of these theorems is that, fixing a base disk $D_0$, one can
start at any disk $D$ and move steadily ``closer'' to $D_0$, in an
appropriate sense, by repeatedly surgering $D$ along an intersection with
$D_0$ that lies outermost on $D_0$. Doing this with a certain amount of
care allows one to produce a null-homotopy for any simplicial map of $S^k$
into the complex. The arguments given in~\cite{JDiff} are basically
correct, but contain some minor misstatements. A much improved treatment is
given in~\cite{Cho}.

\section{$\D(H)$ and the tree $\widetilde{T}$}

Fix a standard unknotted genus-$2$ handlebody $H\subset S^3$. From now on,
disks in $H$ are assumed to be \textit{nonseparating} unless explicitly
stated otherwise. Since at most three nonseparating isotopy classes in $H$
may be represented by disjoint disks, $\D(H)$ has dimension~$2$.

\begin{figure}
\begin{center}
\includegraphics[width=28 ex]{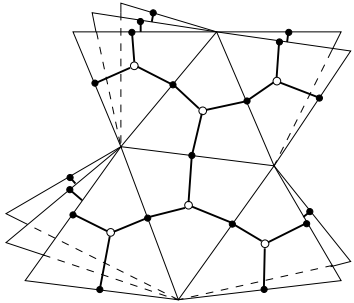}
\caption{A portion of the nonseparating disk complex $\D(H)$ and the tree
$\widetilde{\T}$. Countably many $2$-simplices meet along each edge.}
\label{fig:subdivision}
\end{center}
\end{figure}
A portion of $\D(H)$ is shown in figure~\ref{fig:subdivision}. An edge of
$\D(H)$ is a pair of disjoint nonisotopic disks in $H$, and is called a
\textit{meridian pair,} or just a \textit{pair.}  Similarly, a $2$-simplex
of $\D(H)$ is a \textit{triple.} A triple corresponds to a
\textit{$\theta$-curve} in $H$; that is, a union of three arcs meeting only
in their common endpoints that is a deformation retract of $H$. It is
determined by the condition that each arc passes through exactly one of the
disks of the triple in one point. That arc is called the arc \textit{dual}
to the disk it meets, and the disk is called the disk \textit{dual} to the
arc.

Also shown in figure~\ref{fig:subdivision} is a tree $\widetilde{T}$ which
is a deformation retract of $\D(H)$. It is constructed as follows.  Let
$\D'(H)$ be the first barycentric subdivision of $\D(H)$. Denote by
$\widetilde{\T}$ the subcomplex of $\D'(H)$ obtained by removing the open
stars of the vertices of $\D(H)$. It is a bipartite graph, with ``white''
vertices of valence $3$ represented by triples and ``black'' vertices of
(countably) infinite valence represented by pairs. The valences reflect the
fact that moving along an edge from a triple to a pair corresponds to
removing one of its three disks, while moving from a pair to a triple
corresponds to adding one of infinitely many possible third disks to a
pair.  The possible disjoint third disks that can be added are called the
``slope disks'' for the pair.

The link of a vertex of $\D(H)$ in $\D'(H)$ is an infinite graph contained
in $\widetilde{\T}$, and structurally very similar to $\widetilde{\T}$,
except that its white vertices have valence $2$ rather than
valence~$3$.

It was proven in \cite[Theorem~5.5]{JDiff} and in~\cite{Cho} that the link
in $\D'(H)$ of each vertex of $\D(H)$ is contractible. Combined with
theorem~\ref{thm:contractible}, this implies that
\begin{theorem} $\widetilde{\T}$ is a tree.
\end{theorem}

It is well known that complexes such as the disk complex admit actions of
mapping class groups, indeed this is one of the important motivations for
studying them. The mapping class group $\pi_0(\Diff_+(H))$ acts on $\D(H)$,
simply by $[f]\cdot \langle [D_0],\ldots,[D_k]\rangle = \langle
[f(D_0)],\ldots,[f(D_k)]\rangle$. It is rather obvious that the quotient of
$\widetilde{\T}$ is a single edge, so as seen in~\cite{JDiff}, the
Bass-Serre theory of group actions on trees shows that $\pi_0(\Diff_+(H))$
is a free product with amalgamation. In particular, it is finitely
generated, and an explicit presentation can be worked out by examination of
the vertex and edge stabilizers.

\section{The Goeritz groups and the Scharlemann-Akbas tree}
\label{sec:Goeritz}

The \textit{Goeritz group} $\mathcal{G}$ is the group of isotopy classes of
orientation-preserv\-ing diffeomorphisms of $S^3$ that leave $H$ invariant.
Since every orientation-preserving diffeomorphism of $S^3$ is isotopic to
the identity, $\G$ is exactly the indeterminacy when one takes an arbitrary
unknotted handlebody in $S^3$ and moves it to $H$ by an isotopy.

At times, mostly in section~\ref{sec:links}, we will use the
\textit{extended Goeritz group}~$\G_\pm$, in which the diffeomorphisms of
$S^3$ are allowed to reverse orientation.

Some very important recent work of M. Scharlemann and E. Akbas provides a
precise description of~$\G$. Their methodology is also of interest, and as
we will see it can be considerably simplified by using~$\D(H)$.

L. Goeritz \cite{Goeritz} gave generators for $\mathcal{G}$, and
M. Scharlemann \cite{ScharlemannTree} provided a modern proof that they do
generate. That proof uses a $2$-dimensional complex whose vertices are
isotopy classes of \textit{splitting spheres} for $H$, that is, $2$-spheres
in $S^3$ that intersect $H$ in one essential disk (necessarily a separating
disk), up to isotopy through such spheres. This complex, which we call the
\textit{Scharlemann complex,} is rather difficult to work with, since the
adjacency condition for two spheres in the complex is not disjointness---
this condition would not work, because vertices with disjoint
representatives are equal, and the complex would be a discrete
set. Instead, one must use minimal intersection, that is, the intersections
of the spheres with $\partial H$ meet in only $4$ points. The Goeritz group
acts on the Scharlemann complex in the usual way, with quotient a finite
complex and with finitely generated stabilizers. Using highly nontrivial
geometric arguments, Scharlemann proved:
\begin{theorem}[M. Scharlemann] The Scharlemann complex is connected.\par
\end{theorem}
\noindent It follows by well-known algebraic considerations that
$\mathcal{G}$ is finitely generated.

By additional complicated geometric arguments, E. Akbas
\cite{Akbas} showed the following:
\begin{theorem}[E. Akbas] The natural $1$-dimensional deformation retract
of the Scharlemann complex is a tree.
\end{theorem}
\noindent Akbas worked out the vertex and edge stabilizers of the tree,
enabling him to give a pleasantly transparent presentation
of~$\mathcal{G}$.
\longpage

It turns out that the work of Scharlemann and Akbas is intimately related
to the disk complex, and begins to give insight into the role of
$\widetilde{\T}$. Recall that a disk $\tau$ in $H$ is called
\textit{primitive} when there exists a properly imbedded disk $\tau'$ in
the complementary handlebody $\overline{S^3-H}$ such that the circles
$\partial\tau$ and $\partial\tau'$ in $\partial H$ intersect transversely
in a single point. We call $\tau'$ a \textit{dual disk} of~$\tau$.

A \textit{primitive pair} is a pair of disjoint (nonisotopic) primitive
disks. \textit{The splitting spheres for $H$ correspond exactly to the
primitive pairs.} For a splitting sphere cuts $H$ into two unknotted solid
tori, each containing a unique nonseparating disk which is primitive. On
the other hand, any primitive pair can be shown \cite[Lemma~2.2]{Cho} to
have a unique pair of dual disks that are disjoint from each other and each
disjoint from the other disk of the primitive pair, and the splitting
sphere is the boundary of a small regular neighborhood of the union of
either of the disks with its dual.

The \textit{primitive subtree} $\widetilde{\mathcal{T}}_0$ is the
subcomplex of $\widetilde{\mathcal{T}}$ spanned by the vertices that are
primitive pairs and primitive triples. It is routine to check that two
splitting spheres represent adjacent vertices in the Scharlemann complex if
and only if the corresponding primitive pairs are contained in a primitive
triple.
\begin{figure}
\begin{center}
\includegraphics[width=36 ex]{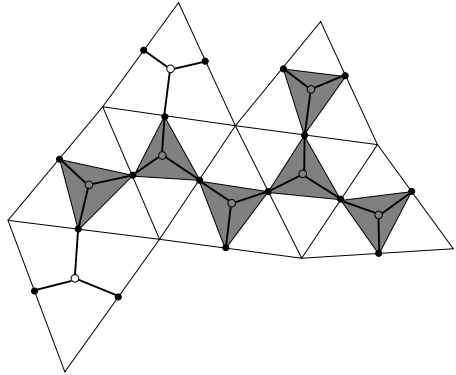}
\caption{A portion of the nonsimplicial imbedding of the Scharlemann
complex into $\D(H)$. The vertices of the complex map to the vertices of
$\D'(H)$ that are primitive pairs.}
\label{fig:primitive_tree}
\end{center}
\end{figure}
So sending a splitting sphere to its corresponding primitive pair
determines a (non-simplicial) imbedding of the Scharlemann complex into
$\D(H)$, as shown in figure~\ref{fig:primitive_tree}, and there is an
obvious deformation from the image to the primitive tree. In this way,
\textit{the Scharlemann-Akbas tree is naturally identified with the
primitive subtree $\widetilde{\T}_0$ of~$\widetilde{\T}$.}

Notice that our observations to this point show that Scharlemann's
connectedness theorem immediately implies Akbas's tree theorem, since any
connected subcomplex of a tree is a tree; also, Scharlemann's theorem
verifies that the primitive subtree really is a tree. But there is now an
independent proof that $\widetilde{\T}_0$ is a tree, based on the following
key fact about primitive disks~\cite[Theorem 2.3]{Cho}:
\begin{theorem} Let $\tau$ and $\sigma$ be primitive disks in $H$ which
intersect transversely. Let $\tau_1$ and $\tau_2$ be the disks that result
from surgering $\tau$ along an intersection arc which is outermost on
$\sigma$. Then $\tau_1$ and $\tau_2$ are primitive.\par
\label{thm:Cho}
\end{theorem}
\noindent The proof of theorem~\ref{thm:Cho} is not trivial, but neither is
it lengthy. A key ingredient is an algebraic fact: if $x$ and $y$ freely
generate $\pi_1(H)$, then a cyclically reduced word that contains both $x$
and $x^{-1}$ cannot be a primitive element (i.~e.~cannot be an element of a
generating pair of $\pi_1(H)$). Theorem~\ref{thm:Cho} easily implies that
$\widetilde{\T}_0$ is contractible, so implies the theorem of Akbas, and it
is straightforward to work out the stabilizers and recover his presentation
of~$\G$.

We remark that theorem~\ref{thm:Cho} is not true for higher genera, so does
not provide a proof that the higher-genus Goeritz groups are finitely
generated.

\section{Tunnels and the tree $\T$}
\label{sec:tunnels}

Consider a knot tunnel, regarded as a $1$-handle attached to a regular
neighborhood of the knot. Since $S^3$ has a unique Heegaard splitting of
each genus, we may move the neighborhood of the knot and the $1$-handle by
isotopy to be the standard $H$. The cocore $2$-disk of the $1$-handle is
then a nonseparating disk in $H$, i.~e.~a vertex of~$\D(H)$. Allowing for
the indeterminacy in this process measured by the Goeritz group, we have
our definition of a knot tunnel:

\begin{definition} A \textit{tunnel} is a $\mathcal{G}$-orbit of disks
in $H$. Thus the tunnels correspond exactly to the vertices of the quotient
$\D(H)/\G$ of $\D(H)$ by the action of $\mathcal{G}$.
\end{definition}

\begin{figure}
\begin{center}
\includegraphics[width=70 ex]{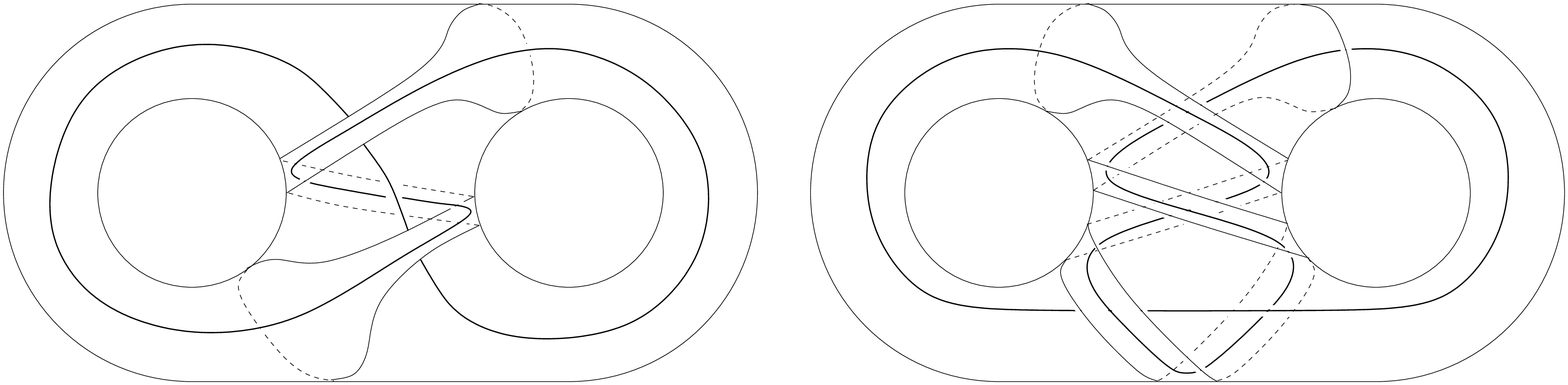}
\caption{Knot tunnels for the right-handed trefoil and figure-8 knots.}
\label{fig:trefoil}
\end{center}
\end{figure}

\begin{notation}
If $\tau$ is a tunnel, then it is a \textit{tunnel of the knot} $K_\tau$
which is a core of the solid torus obtained by cutting $H$
along~$\tau$. We regard $K_\tau$ as defined only up to isotopy in~$S^3$.
\end{notation}
\noindent Figure~\ref{fig:trefoil} shows tunnels for the right-handed
trefoil and figure-8 knots.

One may also consider tunnels and knots up to equivalence that may be
orientation-reversing, simply by replacing $\G$ by~$\G_\pm$ in the
definition of tunnel. This requires only very minor modifications to the
theory. In our exposition, we usually consider only orientation-preserving
equivalence, but will occasionally point out what happens in the other
case.

Our definition of tunnel agrees with the standard definition of
``equivalent'' knot tunnels using a $1$-handle attached to a regular
neighborhood of a knot, and hence with any definition, but it is still
worth thinking through the tunnel arc viewpoint.  Given a tunnel arc
attached to the knot, take a regular neighborhood of the knot and arc
and move it to $H$ by isotopy. The knot and tunnel arc form a
$\theta$-curve, and the disks dual to the arcs of the $\theta$-curve form a
triple. Conversely, any triple containing $\tau$ determines such a
$\theta$-curve whose arcs not dual to $\tau$ form (a knot isotopic to)
$K_\tau$. Thus the isotopy classes in $S^3$ of arcs that determine the
tunnel correspond exactly to the white vertices of the link of $\tau$ in
$\D'(H)/\G$. The moves usually called ``sliding'' change the
$\theta$-curve, and correspond to moving through the link of $\tau$
in~$\D'(H)/\G$.

It is not difficult to see that there is a unique $\G$-orbit $\pi_0$ of
primitive disks. Using \cite[Lemma~2.2]{Cho}, one can show that there
is a unique $\G$-orbit $\mu_0$ of primitive disks. Also from~\cite{Cho},
there is a unique orbit $\theta_0$ of primitive triples.

\begin{figure}\begin{center}
\includegraphics[width=28 ex]{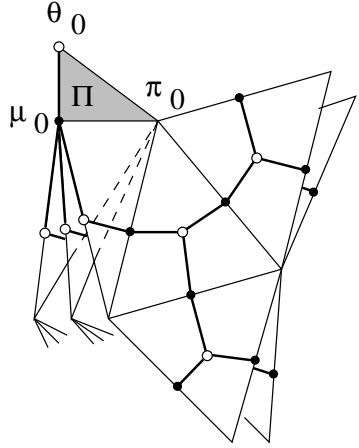}
\caption{A portion of $\D(H)/\G$ and $\T$ near the primitive orbits.}
\label{fig:Delta}
\end{center}
\end{figure}

The quotient $\D'(H)/\G$ inherits a triangulation from~$\D'(H)$. The
quotient of $\widetilde{\mathcal{T}}$ by the action of $\mathcal{G}$ is
a subcomplex of $\D'(H)/\G$. We call this quotient $\T$. A portion of
$\D(H)/\G$ and $\T$ is shown in figure~\ref{fig:Delta}.
\begin{theorem} $\T$ is a tree.
\label{thm:T_is_a_tree}
\end{theorem}
\begin{proof}
Since all primitive triples are equivalent, and the three pairs in a triple
are also equivalent under the stabilzer of the triple, the quotient of
$\widetilde{\T}_0$ by $\G$ is the single edge $\langle
\mu_0,\theta_0\rangle$. Each point $p$ of $\widetilde{\T}-\widetilde{\T}_0$
has a well-defined minimal distance from $\widetilde{\T}_0$, which is just
the length of the unique shortest arc in $\widetilde{\T}$ from $p$ to
$\widetilde{\T}_0$ (as is usual, we use a path metric with the length of
each edge equal to $1$). This arc must end at a primitive pair, since all
edges incident to a primitive triple lie in $\widetilde{\T}_0$. The action
of $\mathcal{G}$ leaves $\widetilde{\T}_0$ invariant, hence preserves the
shortest distance, so the image of the arc from $p$ to $\widetilde{\T}_0$
is the unique shortest arc from the image of $p$ to
$\widetilde{\T}_0$. Thus $\T$ has a unique arc from each point to $\mu_0$,
so $\T$ is a tree.
\end{proof}

In section~\ref{sec:Akbas} we will analyze $\D(H)/\G$ and $\T$ in quite a
bit more detail, after setting up some useful terminology in the next two
sections.

\section{Slope disks, cables, and waves}
\label{sec:slope_disks}

\begin{figure}
\begin{center}
\includegraphics[width=68 ex]{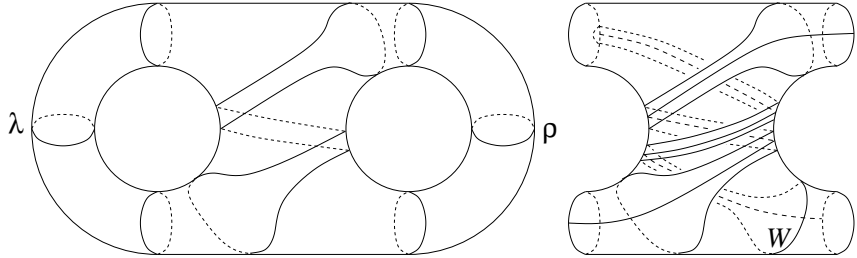}
\caption{A slope disk and associated structures. $W$ is one of the four
waves that correspond to this disk.}
\label{fig:slopes}
\end{center}
\end{figure}

Throughout this section, we consider a pair of disks $\lambda$ and $\rho$
(for ``left'' and ``right'') in $H$, as shown
abstractly in figure~\ref{fig:slopes}. Since $\lambda$ and $\rho$ are
arbitrary, the true picture in $S^3$ might look very different from the
standard-looking pair shown here. Let $B$ be $H$ cut along $\lambda\cup
\rho$. The frontier of $B$ in $H$ consists of four disks which appear
vertical in figure~\ref{fig:slopes}. Denote this frontier by $F$, and let
$\Sigma$ be $B\cap \partial H$, a sphere with four holes.

\begin{definition}
A \textit{slope disk for $\{\lambda,\rho\}$} is an essential disk, possibly
separating, which is contained in $B$ and not isotopic to any component
of~$F$.
\end{definition}
\noindent The boundary of a slope disk always separates $\Sigma$ into two
pairs of pants, conversely any loop in $\Sigma$ that is not homotopic into
$\partial \Sigma$ is the boundary of a unique slope disk. If two slope
disks are isotopic in $H$, then they are isotopic in~$B$.

An arc in $\Sigma$ whose endpoints lie in two different boundary circles of
$\Sigma$ is called a \textit{cabling arc,} and a pair of disjoint cabling
arcs whose four endpoints lie in the four different boundary circles of
$\Sigma$ is called a \textit{cable.}  Figure~\ref{fig:slopes} shows a cable
disjoint from a slope disk. A slope disk is disjoint from a unique
cable. On the other hand, each cabling arc $\alpha$ determines a unique
slope disk: if the endpoints of $\alpha$ lie in the frontier disks $F_1$
and $F_2$ of $B$, then the frontier of a regular neighborhood of $F_1\cup
\alpha\cup F_2$ in $B$ is the slope disk. Finally, each cabling arc
determines a unique cable, the cable disjoint from the slope disk that it
determines.

A \textit{wave} is a disk in $B$ that meets $F$ in a single arc and is
essential in $(B,F)$, that is, not parallel through such disks to a disk in
$F$. A wave $W$ determines a unique slope disk: if $F_0$ is the component
of $F$ that meets $W$, then the frontier of a regular neighborhood of
$F_0\cup W$ in $B$ consists of two disks, one a slope disk and the other
parallel to a component of $F$. A cabling arc $\alpha$ determines two
waves: if $F_1$ and $F_2$ are the disks of $F$ that contain a boundary
point of $\alpha$, the frontier of a regular neighborhood of each
$\alpha\cup F_i$ in $B$ is a wave which determines the same slope disk as
$\alpha$ does. A slope disk is produced by any of four waves, the two pairs
produced from a cabling that determines the slope disk. One such wave is
shown in figure~\ref{fig:slopes}.

In summary:
\begin{enumerate}
\item A slope disk determines a cable, either of whose cabling arcs
determines the slope disk and hence the other cabling arc of the cable.
\item A slope disk determines four waves, each of which determines the
slope disk.
\end{enumerate}

Since disjoint slope disks in $B$ are parallel, disjoint waves in $B$
determine the same slope disk.

\section{General slope coordinates}
\label{sec:general_slopes}

In this section, we will explain how each choice of nonseparating slope
disk for a pair $\mu=\{\lambda,\rho\}$ determines a correspondence
between~$\Q\cup\{\infty\}$ and the set of all slope disks of $\mu$. As we
saw in section~\ref{sec:slope_disks}, such a correspondence associates a
value to each cabling arc, cable, and wave as well. In fact, it will be
most convenient to use cabling arcs to define the value, as the method is
simply a version of the standard procedure for associating a parameter to a
rational tangle.

\begin{figure}
\begin{center}
\includegraphics[width=45 ex]{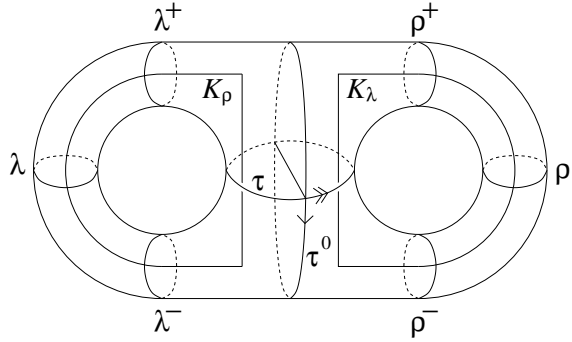}
\caption{The slope-zero perpendicular disk $\tau^0$. It is chosen so that
$K_\lambda$ and $K_\rho$ have linking number~$0$.}
\label{fig:slope_coords}
\end{center}
\end{figure}

Let $F$ and $\Sigma$ be as in the previous section. Fixing a nonseparating
slope disk $\tau$ for $\mu$, we will write $(\mu;\tau)$ for the ordered
pair consisting of $\mu$ and $\tau$.
\begin{definition} A \textit{perpendicular disk} for $(\mu;\tau)$
is a disk $\tau^\perp$, with the following properties:
\begin{enumerate}
\item $\tau^\perp$ is a slope disk for $\mu$.
\item $\tau$ and $\tau^\perp$ intersect transversely in one arc.
\item $\tau^\perp$ separates $H$.
\end{enumerate}
\end{definition}
There are infinitely many choices for $\tau^\perp$, but because $H\subset
S^3$ there is a natural way to choose a particular one, which we call
$\tau^0$. It is illustrated in figure~\ref{fig:slope_coords}. To construct
it, start with any perpendicular disk and change it by Dehn twists of $H$
about $\tau$ until the core circles of the complementary solid tori have
linking number~$0$.

For calculations, it is convenient to draw the picture as in
figure~\ref{fig:slope_coords}, and orient the boundaries of $\tau$ and
$\tau^0$ so that the orientation of $\tau^0$ (the ``$x$-axis''), followed
by the orientation of $\tau$ (the ``$y$-axis''), followed by the outward
normal of $H$, is a right-hand orientation of $S^3$. At the other
intersection point, these give the left-hand orientation, but we will see
that the coordinates are unaffected by changing the choices of which of
$\{\lambda,\rho\}$ is $\lambda$ and which is $\rho$, or changing which
sides are $+$ and which are $-$, provided that the $+$ sides both lie on
the same side of $\lambda\cup\rho\cup\tau$ in
figure~\ref{fig:slope_coords}.

\begin{figure}
\begin{center}
\includegraphics[width=\textwidth]{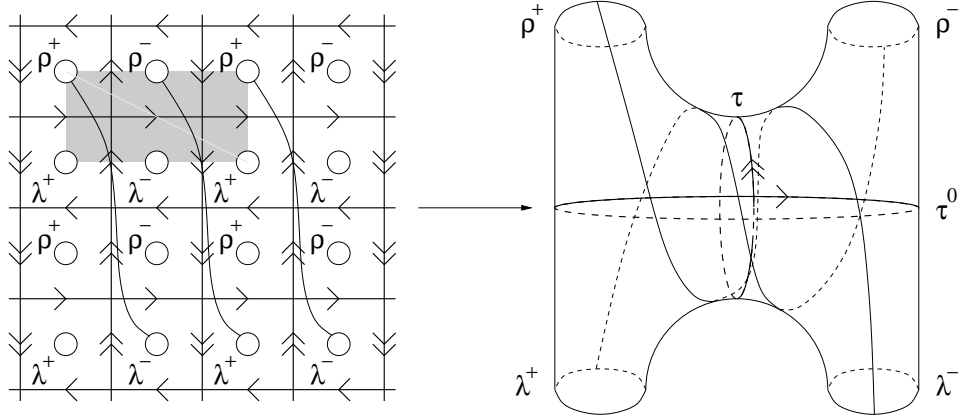}
\caption{The covering space $\widetilde{\Sigma}\to\Sigma$, and some lifts
of a $[1,-3]$-cabling arc. The shaded region is a fundamental domain.}
\label{fig:covering}
\end{center}
\end{figure}

Let $\widetilde{\Sigma}$ be the covering space of $\Sigma$ such that:
\begin{enumerate}
\item $\widetilde{\Sigma}$ is the plane with an open disk of radius $1/8$
removed from each point with half-integer coordinates.
\item The components of the preimage of $\tau$ are the vertical lines
with integer $x$-coordinate.
\item The components of the preimage of $\tau^0$ are the horizontal lines
with integer $y$-coordinate.
\end{enumerate}
\noindent Figure~\ref{fig:covering} shows a picture of $\widetilde{\Sigma}$
and a fundamental domain for the action of its group of covering
transformations, which is the orientation-preserving subgroup of the group
generated by reflections in the half-integer lattice lines (that pass
through the centers of the missing disks). Each circle of
$\partial\widetilde{\Sigma}$ double covers a circle of~$\partial \Sigma$.

If we lift any cabling arc in $\Sigma$ to $\widetilde{\Sigma}$, the lift
runs from a boundary circle of $\widetilde{\Sigma}$ to one of its
translates by a vector $(p,q)$ of signed integers, defined up to
multiplication by the scalar $-1$. Thus each cabling arc receives a
\textit{slope pair} $[p,q]=\{(p,q),(-p,-q)\}$, and is called a
\textit{$[p,q]$-cabling arc.} Of course, each slope disk, cable, or wave
receives a slope pair $[p,q]$ as well.

An important observation is that a $[p,q]$-slope disk is nonseparating in
$H$ if and only if $q$ is odd. Both happen exactly when the corrsponding
cabling arc has one endpoint in $\lambda^+$ or $\lambda^-$ and the other in
$\rho^+$ or~$\rho^-$.

\begin{definition} The \textit{$(\mu;\tau)$-slope} of a $[p,q]$-slope 
disk, cabling arc, cable, or wave is the rational number $q/p$.
\end{definition}

Finally, we clarify why the choices of $\lambda$ and $\rho$ and of $+$ and
$-$ sides do not affect the final slope. Interchanging $\lambda$ and $\rho$
changes the covering space in figure~\ref{fig:covering} by a vertical
translation by $1$, while interchanging the $+$ and $-$ sides changes it by
a horizontal translation by $1$. The orientation on one of $\tau^0$ or
$\tau$ is reversed, but the set of lifts in figure~\ref{fig:covering} is
preserved, so the slope pairs are unchanged.

\section{Slope disks of primitive pairs}
\label{sec:Akbas}

In this section, we will examine the set of slope disks for primitive
pairs. In particular, we will show that the equivalence classes modulo $\G$
correspond, by taking the ``simple slope,'' to $\Q/\Z\cup\{\infty\}$.

Fix a primitive pair $\mu_0=\{\lambda_0,\rho_0\}$. We use the notation of
the previous section, but add ``$0$'' subscripts as in $\Sigma_0$ and
$\widetilde{\Sigma}_0$ to remind ourselves that we are in the primitive
case.

First, we review Akbas's \cite{Akbas} description of the stabilizer
$\G_{\mu_0}$ of $\mu_0$ under the action of $\G$. Three diffeomorphisms of
$H$ are the ``hyperelliptic'' involution $\alpha$ that fixes an arc in each
of $\{\lambda_0,\rho_0\}$ in $H$, the left-hand ``half-twist'' $\beta$ that
fixes $\lambda_0$ and reflects $\rho_0$ across an arc, and a ``rotation''
involution $\gamma$ that interchanges $\lambda_0$ and $\rho_0$ and fixes an
arc from the reader's nose straight through the middle of $H$.  As usual,
we make no notational distinction between these maps and their isotopy
classes.

As is well known, the restriction of $\alpha$ to $\partial H$ is central in
the mapping class group of $\partial H$ and acts trivially on every simple
closed loop in~$\partial H$. Consequently $\alpha$ is a central involution
in $\G$ which acts trivially on every disk in~$H$. It is evident that
$\beta\gamma\beta^{-1}=\alpha\gamma$. From~\cite{Akbas} and also, using the
present viewpoint, from~\cite{Cho}, we have:
\begin{proposition} The stabilizer $\mathcal{G}_{\mu_0}$ is the subgroup
generated by $\alpha$, $\beta$, and $\gamma$.  In fact, $\G_{\mu_0}$ is
the semidirect product $(C_2\times \Z)\circ C_2$, where $\langle
\alpha,\beta\rangle$ is the normal subgroup $C_2\times \Z$ and $\gamma$
acts by $\gamma\alpha\gamma^{-1}=\alpha$ and
$\gamma\beta\gamma^{-1}=\alpha\beta$.
\label{prop:stabilizer}
\end{proposition}

Now we analyze the action of $\G$ on the slope disks associated to
$\mu_0$. Each cable for $\mu_0$, and hence each slope disk for $\mu_0$, is
invariant under $\alpha$ and $\gamma$. In terms of the covering space
$\widetilde{\Sigma}_0$ from the previous section, $\alpha$ lifts to a
horizontal translation by $1$, and $\gamma$ lifts to multiplication by
$-I$.

To understand the action of $\beta$, fix some primitive slope disk $\pi_0$
for $\mu_0$, and consider $(\mu_0;\pi_0)$-slope pairs.  The associated
perpendicular disk $\pi_0^0$ is exactly the intersection of $H$ with the
splitting sphere disjoint from $\lambda_0\cup \rho_0$. A lift of $\beta$ to
$\widetilde{\Sigma}_0$ sends $(x,y)$ to $(x+y,y)$, so $\beta$ sends a
cabling arc with slope pair $[p,q]$ to one with slope pair~$[p+q,q]$.

The next observation is that the action of $\G_{\mu_0}$ is transitive on
primitive slope disks for $\mu_0$. Perhaps the easiest way to see this is
that the $\beta$-translates of $\pi_0$ are the only slope disks for $\mu_0$
whose complementary tori in $H$ are unknotted in $S^3$ (for the others, the
core circle of the complementary torus is a nontrivial $2$-bridge knot, see
section~\ref{sec:simple_tunnels}).

We claim that slope disks of $\mu_0$ are equivalent under $\G$-orbits only
when they are equivalent under $\G_{\mu_0}$. We have just seen the action
of $\G_{\mu_0}$ on the primitive slope disks of $\mu_0$ is transitive.
Consider two nonprimitive slope disks $\tau_1$ and $\tau_2$ of $\mu_0$ that
are equivalent by an element $g$ of $\G$. The shortest path $\omega_i$ in
$\widetilde{\T}$ from any vertex in the link of $\tau_i$ in $\D'(H)$ to a
vertex of the primitive tree $\widetilde{\T}_0$ is the edge of
$\widetilde{\T}$ from $\mu_0\cup\{\tau_i\}$ to $\mu_0$. Since the primitive
tree is invariant, $g$ must take $\omega_1$ to $\omega_2$ and hence $g\in
\G_{\mu_0}$.
\longpage

Since $\alpha$ and $\gamma$ act trivially on slope disks associated to
$\mu_0$, while $\beta$ sends a $[p,q]$-disk to a $[p+q,q]$-disk, it follows
that \textit{sending a $[p,q]$-slope disk associated to $\mu_0$ to $p/q$
induces a bijection from the set of $\G$-orbits of slope disks associated
to $\mu_0$ to the set $\Q/\Z\cup\{\infty\}$.} The slope disk that
corresponds to $\infty$ is the one with slope pair $[1,0]$, that
is,~$\pi_0^0$.

Since any other choice of $\pi_0$ differs from our previous one by the
action of a power of $\beta$, the bijection to $\Q/\Z\cup\{\infty\}$ is
independent of the choice of $\pi_0$. Moreover, this shows that the
primitive slope disks are those having slope pair of the form $[p,1]$, so
they correspond to $[0]\in\Q/\Z\cup\{\infty\}$.

The slope disks for a primitive pair play a key role in our theory, so we
introduce some special terminology.
\begin{definition} A possibly separating disk in $H$ is called \textit{simple} 
if it is nonprimitive and is a slope disk for a primitive pair.
\end{definition}

As we have seen, each simple disk is preserved up to isotopy by
$\alpha$ and $\gamma$, while $\beta$ sends a simple disk of slope pair
$[p,q]$ to one with slope pair $[p+q,q]$.

\begin{definition}
The \textit{simple slope} of a simple disk is the corresponding element
$[p/q]\in (\Q/\Z-\{0\})\cup\{\infty\}$, which has $q$~odd if and only if
the disk is nonseparating.
\end{definition}

\section{The tree of knot tunnels}
\label{sec:quotients}
\shortpage

In this section we will analyze the quotient $\D(H)/\G$ and the tree
$\T$. We begin with a summary description, which we will establish in the
remainder of this section. It may be useful at this point to refer to
figure~\ref{fig:Delta}. To fix notation, let $\theta_0$ be a primitive
triple $\{\lambda_0,\rho_0,\pi_0\}$, containing the primitive pair
$\mu_0=\{\lambda_0,\rho_0\}$.
\begin{enumerate}
\item[(1)] The quotient of the subcomplex of $\D(H)$ spanned by primitive
vertices is a single simplex $\Pi$ spanned by $\pi_0$, $\mu_0$, and
$\theta_0$. It is a copy of a $2$-simplex of $\D'(H)$, and meets $\T$ in
the edge $\langle\mu_0,\theta_0\rangle=\widetilde{\T}_0/\G$.
\item[(2)]
Each $2$-simplex of $\D(H)$ having two primitive vertices and one simple
vertex gets identified with some other such simplices, then folded in half
by the action of some conjugate of the element $\gamma$ of
section~\ref{sec:Akbas}, and the quotient half-simplex is attached to $\Pi$
along the edge $\langle \pi_0,\mu_0\rangle$. These half-simplices
correspond to the elements of $\Q/\Z-\{0\}$ with odd denominator. The
intersection of $\T$ with such a half-simplex is a pair of edges.
\item[(3)]
Each component of the remainder of $\D(H)$ descends injectively to
$\D(H)/\G$. In fact, it is identified with some other such components, and
the result is attached to a half-simplex of $\D(H)/\G$ along the edge from
$\pi_0$ to the simple vertex.  The intersection of $\T$ with each component
of the remainder of $\D(H)/\G$ is exactly a copy of the intersection of
$\widetilde{\T}$ with a component of the remainder of~$\D(H)$.
\end{enumerate}

Now we present more detail of this description. Part~(1) should be clear
since the action of $\G$ is transitive on primitive triples, and the
stabilizer of $\theta_0$ acts as the permutation group on the three disks
in~$\theta_0$.

Consider a $2$-simplex of $\D(H)$ which has $\lambda_0$ and $\rho_0$ as two
of its vertices, and third vertex a simple disk. We saw in
section~\ref{sec:Akbas} that the generator $\alpha$ of $\G_{\mu_0}$ acts
trivially on every $2$-simplex, and the generator $\beta$ preserves each of
$\lambda_0$ and $\rho_0$, but acts on slope disks for $\mu_0$ by sending
the disk with slope pair $[p,q]$ to the one with slope pair $[p+q,q]$.
Finally, $\gamma$ interchanges $\lambda_0$ and $\rho_0$, while preserving
each slope disk associated to $\mu_0$. So the effect of $\G_{\mu_0}$ on the
collection of $2$-simplices of $\D(H)$ having $\lambda_0$, $\rho_0$, and a
simple disk as vertices is to identify those for which the ratio $p/q$
differs by an integer, and to fold the resulting $2$-simplices in half,
producing half-sized simplices in $\D(H)/\G$ with short edges attached to
$\Pi$ along the $1$-simplex $\langle \pi_0, \mu_0\rangle$. There is one
such half-simplex for each orbit of nonseparating simple disks under
$\G_{\mu_0}$, and as seen in section~\ref{sec:Akbas}, these orbits
correspond to the elements of $\Q/\Z -\{0\}$ with odd denominator.

To understand the remaining portions of $\D(H)/\G$, fix a simple disk
$\tau$, and consider a portion of $\D(H)$ attached to the $2$-simplex
$\langle \lambda_0,\rho_0,\tau\rangle$ along the edge
$\langle\rho_0,\tau\rangle$. We know that the action of $\beta$ moves it to
similar portions attached along edges of other $2$-simplices having
$\lambda_0$ and $\rho_0$ as vertices. The action of $\gamma$ interchanges
$\lambda_0$ and $\rho_0$, so interchanges the portion with another one
attached along $\langle\lambda_0,\tau\rangle$. Finally, $\alpha$ acts
trivially. So the original portion descends injectively into $\D(H)/\G$
onto a copy attached along the edge~$\langle \pi_0, \tau\rangle$.

The quotient of the primitive subtree $\widetilde{\T}_0$ is the edge
$\langle \theta_0,\mu_0\rangle$ of $\Pi$. The rest of $\T$ should be clear
from the previous discussion and figure~\ref{fig:Delta}. In $\T$, the
vertex $\theta_0$ has valence~$1$, and $\mu_0$ has countable valence. Each
additional edge emanating from $\mu_0$ goes to a white vertex of
valence~$2$, representing a triple with two primitive and one simple disk
as vertices. From then on, $\T$ looks exactly like a portion
of~$\widetilde{\T}$.

We now have a complete picture of $\D(H)/\G$ and $\T$. It will also be
useful to know the stabilizers of the action of $\G$ on nonprimitive disks
in $H$.
\begin{proposition} The stabilizer in $\G$
of a nonseparating, nonprimitive disk $\tau$ in $H$ is as follows:
\begin{enumerate}
\item[(i)] If $\tau$ simple, then its stabilizer is conjugate to the
subgroup $C_2\times C_2$ generated by $\alpha$ and~$\gamma$.
\item[(ii)] Otherwise, the stabilizer is $C_2$, generated by $\alpha$.
\end{enumerate}
\label{prop:tau_stabilizer}
\end{proposition}

\begin{proof} Suppose first that $\tau$ simple. There is a unique
primitive pair $\mu$ which together with $\tau$ spans a $2$-simplex of
$\D(H)$, for if there were two such pairs $\mu$ and $\mu'$, then the
$1$-simplices $\langle\mu\cup \{\tau\},\mu\rangle$ and 
$\langle \mu'\cup \{\tau\},\mu'\rangle$, connected by an
arc in the primitive tree, would form an imbedded loop in
$\widetilde{\T}$. By conjugation, we may assume that $\mu$ is the
standard primitive pair~$\mu_0$.

Let $h$ be in $\G_\tau$. Since $h$ also preserves the primitive
subcomplex of $\D(H)$, it must preserve $\mu_0$. From
proposition~\ref{prop:stabilizer}, any element of $\G_{\mu_0}$ can be
written in the form $\alpha^{\epsilon_1}\beta^n\gamma^{\epsilon_2}$, where
the $\epsilon_i\in\{0,1\}$ and $n\in \Z$. Since $h$ preserves the slope
pair of $\tau$, we have $n=0$ and part~(i) follows.

For~(ii), suppose that $\tau$ is not simple. Let $\theta$ be the vertex of
the link of $\tau$ in $\D'(H)$ that is closest in $\widetilde{\T}$ to the
primitive subtree $\widetilde{\T}_0$. Since the action of $\G$ fixes $\tau$
and preserves $\widetilde{\T}_0$, it must fix the path from $\theta$ to the
nearest point $\mu$ of $\widetilde{\T}_0$ (we will see in
section~\ref{sec:parameterization} that the image of this path in
$\D(H)/\G$ is the \textit{principal path} of $\tau$, illustrated in
figure~\ref{fig:principal_vertex}). By conjugation, we may assume that
$\mu=\mu_0$. Write this path as $\mu_0$, $\mu_0\cup\{\tau_0\}$, $\mu_1$,
$\mu_1\cup \{\tau_1\},\ldots\,$, $\theta$. As in part~(i), any element of
$\G_\tau$ fixes this path, so must be in the stabilzier $\langle \alpha,
\gamma\rangle$ of $\tau_0$. But $\gamma$ interchanges the disks of
$\mu_0$, so acts as reflection on the $2$-simplex spanned by $\mu_0\cup
\{\tau_0\}$ and cannot fix $\mu_1$. Part~(ii) follows.
\end{proof}

\begin{corollary} Let $\langle \tau_1,\tau_2,\tau_3\rangle$ be a
$2$-simplex of $\D(H)$, with $\tau_3$ nonprimitive. Suppose that $\langle
\tau_1,\tau_2,\tau_3\rangle$ is stabilized by an element $h$ of $\G$ other
than the identity or $\alpha$. Then $\tau_3$ is simple, $\tau_1$ and
$\tau_2$ are primitive, and $h$ is conjugate to $\gamma$ or $\alpha\gamma$.
\label{coro:simplex_stabilizers}
\end{corollary}

\begin{proof} Of the three vertices $\langle \tau_1,\tau_2\rangle$, $\langle
\tau_1,\tau_3\rangle$, and $\langle \tau_2,\tau_3\rangle$, let $\langle
\tau_i,\tau_j\rangle$ be the one closest to $\widetilde{\T}_0$ (possibly in
$\widetilde{\T}_0$), and let $\tau_k$ be the vertex of $\langle
\tau_1,\tau_2,\tau_3\rangle$ different from $\tau_i$ and $\tau_j$.  The
action of $h$ must preserve $\langle \tau_i,\tau_j\rangle$, so must
fix~$\tau_k$.

Since $\langle\tau_i, \tau_j\rangle$ is the closest vertex to
$\widetilde{T}_0$, $\tau_k$ is nonprimitive.
Proposition~\ref{prop:tau_stabilizer} then shows that $\tau_k$ is simple
and $h$ is conjugate to either $\gamma$ or $\alpha\gamma$. Again since
$\langle\tau_i, \tau_j\rangle$ is the closest vertex to $\widetilde{T}_0$,
$\tau_i$ and $\tau_j$ must both be primitive, so $\tau_k=\tau_3$ and the
corollary is established.
\end{proof}

\section{Simple, semisimple, and regular tunnels}
\label{sec:simple_tunnels}

A $\G$-orbit of nonseparating simple disks is called a \textit{simple
tunnel.}
\begin{proposition} The simple tunnels are exactly the
``upper'' and ``lower'' tunnels of $2$-bridge knots.
\label{prop:top_and_bottom}
\end{proposition}
\begin{proof}
Let $\Sigma_0$ be as in section~\ref{sec:Akbas}. If $\sigma$ is a simple
tunnel, then $K_\sigma$ is isotopic to the union of the cable in $\Sigma_0$
determined by $\sigma$, plus two arcs in $\overline{\partial H-\Sigma_0}$,
each crossing $\lambda_0$ or $\rho_0$ in one point. A dual arc to $\sigma$
is an arc cutting once across $\sigma$ and connecting the cable arcs. This
is a standard description of the upper and lower tunnels of $2$-bridge
knots.
\end{proof}

From sections~\ref{sec:Akbas} and~\ref{sec:quotients}, we have immediately:
\begin{proposition} The simple tunnels are classified up to equivalence 
by their simple slopes in $\Q/\Z-\{0\}$, and up to possibly
orientation-reversing equivalence by the pairs $\{[p/q],[-p/q]\}$ with
$q$~odd.
\label{prop:simple_classification}
\end{proposition}
\noindent Of course, $[p/q]$ is a version of the standard rational
invariant that classifies $2$-bridge knots. This will be examined further
in section~\ref{sec:2bridge}.

Recall that a $\theta$-curve in $S^3$ is called \textit{unknotted} if its
complement is an open handlebody, and \textit{planar} if it is isotopic
into a standard plane. 
The results of section~\ref{sec:quotients} allow us to describe the
automorphisms of nonplanar unknotted $\theta$-curves.
\begin{corollary} Let $\Theta$ be a nonplanar unknotted $\theta$-curve in
$S^3$, and let $h$ be a orientation-preserving homeomorphism of $S^3$ that
preserves $\Theta$. Then, either
\begin{enumerate}
\item[(i)] $h$ is isotopic preserving $\Theta$ to a homeomorphism which is the
identity or the hyperelliptic involution on a neighborhood of $\Theta$, or
\item[(ii)] $\Theta$ is the union of a $2$-bridge knot and one of its simple
tunnel arcs, and $h$ preserves the tunnel arc and interchanges the two
arcs of the knot.
\end{enumerate}
\label{coro:theta_homeos}
\end{corollary}

\begin{proof}
By isotopy we may assume that $\Theta\subset H$ and $h\in \G$. In $\D(H)$,
the dual disks of the arcs of $\Theta$ span a $2$-simplex, which is
preserved by $h$. Since $\Theta$ is nonplanar, the dual disk of the tunnel
is not primitive. Applying corollary~\ref{coro:simplex_stabilizers},
either case~(i) holds, or the tunnel disk is simple and $h$ is conjugate to
$\gamma$ or to $\alpha\gamma$, giving case~(ii).
\end{proof}

\begin{remark} Corollary~\ref{coro:theta_homeos} holds as stated without
the assumption that $h$ is orientation-preserving, by using
corollary~\ref{coro:extended_simplex_stabilizers} below in place of
corollary~\ref{coro:simplex_stabilizers} in the proof. That is, nonplanar
unknotted $\theta$-curves have no orientation-reversing automorphisms.
\end{remark}

Corollary~\ref{coro:theta_homeos} gives immediately a result of D. Futer
\cite{Futer}:
\begin{theorem}[D. Futer] Let $A$ be a tunnel arc for a nontrivial knot
$K\subset S^3$. Then $A$ is fixed pointwise by a strong inversion of $K$ if
and only if $K$ is a two-bridge knot and $A$ is its upper or lower tunnel.
\label{thm:Futer}
\end{theorem}

\begin{definition} A tunnel $\tau$ is called \textit{semisimple} if $\tau$
is not primitive or simple but $\tau$ lies in the link in $\D(H)$ of a
primitive disk. It is called \textit{regular} if it is neither primitive,
simple, or semisimple.
\end{definition}

Recall that a $(1,1)$-knot is a knot of the form $t_1\cup t_2$, where $t_1$
and $t_2$ are trivial (boundary parallel) arcs in the tori $W_1$ and $W_2$
of a genus-$1$ Heegaard splitting of $S^3$. A \textit{$(1,1)$-tunnel} is a
tunnel of a $(1,1)$-knot which is representable by an arc $t$ in $W_1$
which meets $t_1$ in its endpoints and together with the arc that $\partial
t$ bounds in $t_1$ forms a core circle of $W_1$. There is also a
$(1,1)$-tunnel obtained by the corresponding construction in $W_2$, which
may be equivalent to the one in $W_1$. The knot may have more
$(1,1)$-tunnels coming from other $(1,1)$-descriptions, or may have a
regular tunnel (as occurs for most torus knots). The following fact is
well-known, see for example \cite[Proposition~1.3]{Morimoto-Sakuma}.
\begin{proposition}
A tunnel is a $(1,1)$-tunnel for a nontrivial $(1,1)$-knot
if and only if it is simple or semisimple.
\label{prop:semisimple_is_(1,1)}
\end{proposition}

\begin{remark} The torus knots are the only examples we know of knots
having both $(1,1)$-tunnels and a regular tunnel. We do not know whether
there exists a knot that has more than one regular tunnel.
\end{remark}

\section{Principal paths, principal vertices, and parameterization}
\label{sec:parameterization}

\begin{figure}
\begin{center}
\includegraphics[width=55ex]{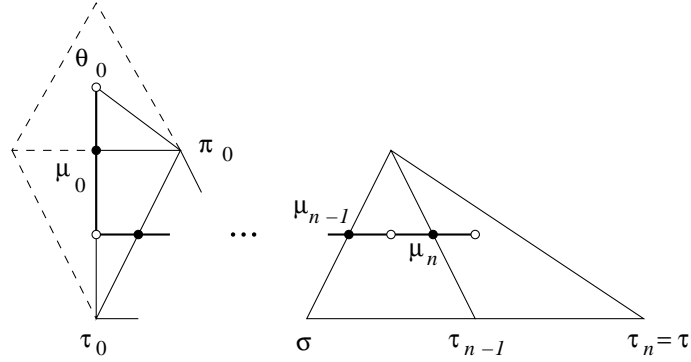}
\caption{The principal path of $\tau$ is the path in $\T/\G$ from
$\theta_0$ to the principal vertex $\mu_n\cup\{\tau\}$ of $\tau$. The
``trailing'' disk $\sigma$, which is the disk of $\mu_{n-1}-\mu_n$,
plays an important role in the calculation of slope invariants.}
\label{fig:principal_vertex}
\end{center}
\end{figure}

Using the coordinates from section~\ref{sec:general_slopes}, we will obtain
a natural numerical parameteri\-zation of all knot tunnels. First, we give
an important definition.
\begin{definition} Let $\tau$ be a tunnel. The
\textit{principal vertex} of $\tau$ is the vertex of the link of $\tau$ in
$\D'(H)/\G$ that is closest to $\theta_0$ in $\T$. The \textit{principal
path} of $\tau$ is the unique path in $\T$ from $\theta_0$ to the principal
vertex of~$\tau$.
\end{definition}
\noindent Figure~\ref{fig:principal_vertex} illustrates the principal path
and principal vertex of a typical tunnel.

\begin{remark}
We will see in lemma~\ref{lem:muplusminus} that the principal vertex of
$\tau$ is $\langle \mu^+,\mu^-,\tau\rangle$, where $\mu^+$ and $\mu^-$ is
the pair of disks used in the definition of the Scharlemann-Thompson
invariant~\cite{Scharlemann-Thompson}. Thus in
figure~\ref{fig:principal_vertex}, $\mu_n=\{\mu^+,\mu^-\}$.
\end{remark}

We can now give the numerical parameterization.
\begin{parameterizationtheorem}
Let $\tau$ be a knot tunnel with principal path $\theta_0$, $\mu_0$,
$\mu_0\cup \{\tau_0\}$, $\mu_1,\ldots\,$, $\mu_n$, $\mu_n\cup \{\tau_n\}$.
Fix a lift of the principal path to $\D(H)$, so that each $\mu_i$
corresponds to an actual pair of disks in~$H$.
\begin{enumerate}
\item If $\tau$ is primitive, put $m_0=[0]\in\Q/\Z$. Otherwise,
let $m_0=[p_0/q_0]\in\Q/\Z$ be the simple slope of $\tau_0$.
\item If $n\geq 1$, then for $1\leq i\leq n$ let $\sigma_i$ be the unique disk
in $\mu_{i-1}-\mu_i$ and let $m_i=q_i/p_i\in\Q$ be the
$(\mu_i;\sigma_i)$-slope of $\tau_i$.
\item If $n\geq 2$, then for $2\leq i\leq n$ define $s_i=0$ or $s_i=1$
according to whether or not the unique disk of $\mu_i\cap\mu_{i-1}$ equals
the unique disk of $\mu_{i-1}\cap\mu_{i-2}$.
\end{enumerate}
Then, sending $\tau$ to the pair $((m_0,\ldots,m_n),(s_2,\ldots,s_n))$ is a
bijection from the set of all tunnels of all tunnel number~$1$ knots to the
set of all elements $(([p_0/q_0],q_1/p_1,\ldots,q_n/p_n),(s_2,\ldots,s_n))$
in
\[\big(\Q/\Z \big)\,\cup\, 
\big(\Q/\Z\,\times\, \Q\big) \,\cup\, \big(\cup_{n\geq 2}
\;\Q/\Z\,\times\, \Q^n \,\times\,\, C_2^{n-1}\big)\]
with all $q_i$ odd.
\label{thm:parameterization}
\end{parameterizationtheorem}
\begin{proof}
Only the sequence $(s_2,\ldots,s_n)$ needs explanation. Referring to
figure~\ref{fig:Delta}, we note that the value of $m_0$ determines a unique
simple tunnel $\tau_0$ (unless $m_0=[0]$, which determines the primitive
tunnel), and the value of $m_1$ determines a unique choice of
$\tau_1$. From then on, one must make a choice of which disk of $\mu_{i-1}$
will be retained in $\mu_i$, and the number $s_i$ records this choice.
\end{proof}

In practice, it is not an easy matter to work out the slope parameters of a
given knot, but in section~\ref{sec:2bridge} we compute them for all
tunnels of $2$-bridge knots, and in~\cite{CM} for the regular tunnels of
torus knots.

\begin{remark}
Theorem~\ref{thm:link_parameterization} below is a version of the
Parameterization Theorem~\ref{thm:parameterization} that includes tunnels
of tunnel number $1$ links. The main difference is that the final slope
$m_n$ may have $q_n$ even.
\end{remark}

\begin{remark}
Reversing the orientation of $S^3$ has the effect of negating each slope
invariant, so the classification of tunnels up to arbitrary homeomorphism
of $S^3$ is obtained from that of the Parameterization
Theorem~\ref{thm:parameterization} by adding the equivalence
$(m_0,m_1,\ldots, m_n)\sim (-m_0,-m_1,\ldots, -m_n)$.
\label{rem:ORequivalence}
\end{remark}

\section{The cabling construction} 
\label{sec:cabling_construction}

In this section, we will see how the tree $\T$ specifies a unique sequence
of ``cabling constructions'' that produce a given tunnel. Roughly speaking, a
path of length $2$ from a white vertex to a white vertex corresponds to one
cabling construction. The principal path of $\tau$ encodes the unique
sequence of cabling constructions that produces~$\tau$.

The cabling construction is simple and will look familiar to experts. In a
sentence, it is ``Think of the union of $K$ and the tunnel arc as a
$\theta$-curve, and rationally tangle the ends of the tunnel arc and one of
the arcs of $K$ in a neighborhood of the other arc of $K$.''  We sometimes
call this ``swap and tangle,'' since one of the arcs in the knot is
exchanged for the tunnel arc, then the ends of other arc of the knot and
the tunnel arc are connected by a rational tangle.

The cabling operation is a very restricted special case of the ``tunnel
moves'' described in Goda, Scharlemann, and Thompson \cite{G-S-T}. A tunnel
move means a replacement of $K_\tau$ by any knot in $\partial H$ that
crosses $\partial\tau$ exactly once. A tunnel move is a composition of an
arbitrarily large number of cabling constructions that retain the same arc
of the knot. As mentioned in section~\ref{sec:overview}, tunnel moves are
studied from our viewpoint in~\cite{CM}.

We begin with some terminology.
\begin{definition} 
Let $\mu$ be a pair, and let $\tau$ be a disk in $H$.  We say that $\mu$ is
a \textit{meridian pair of $\tau$} when $\tau$ is a slope disk of $\mu$. In
this case, the pair $(\mu;\tau)$ corresponds to the directed $1$-simplex in
$\T$ (or $\widetilde{\T}$) from $\mu$ to $\mu\cup \{\tau\}$.
\end{definition}
\noindent The meridian pairs of $\tau$ correspond exactly to the white
vertices of the link of $\tau$ in~$\D'(H)$.

Geometrically, $(\mu;\tau)$ corresponds to an isotopy class of tunnel arc
of $K_\tau$. In a $\theta$-curve corresponding to $\mu\cup \{\tau\}$, the
union of the arcs dual to $\mu$ is $K_\tau$, and the arc dual to $\tau$ is
a tunnel arc for $K_\tau$. Here, isotopy class refers to the isotopy class
in $H$ when working in $\widetilde{\T}$, and the isotopy class in $S^3$
(possibly moving $K_\tau$ along with the arc) when working in~$\T$.

\begin{figure}
\begin{center}
\includegraphics[width=\textwidth]{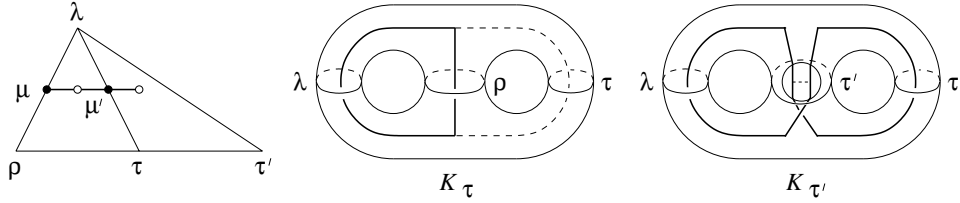}
\caption{Schematic for the general cabling construction. In the middle ball
in the right-hand picture of $H$, the two vertical arcs form some rational
tangle, disjoint from the disk~$\tau'$.}
\label{fig:schematic}
\end{center}
\end{figure}
Moving through the tree determines a sequence of steps in which one of the
two disks of a pair $\{\lambda,\rho\}$ is replaced by a tunnel disk $\tau$,
and a slope disk $\tau'$ of the new pair $\mu'$ (with $\tau'$ nonseparating
in $H$) is chosen as the new tunnel disk. As illustrated in
figure~\ref{fig:schematic}, the way the path determines the particular
cabling operation is:
\begin{enumerate}
\item The selection of $\lambda$ or $\rho$ corresponds to which edge one
chooses to move out of the white vertex $\{\lambda,\rho,\tau\}$.
\item The selection of the new slope disk $\tau'$ corresponds to which edge
one chooses to continue out of the black vertex $\mu'$.
\end{enumerate}

\begin{figure}
\begin{center}
\includegraphics[width=\textwidth]{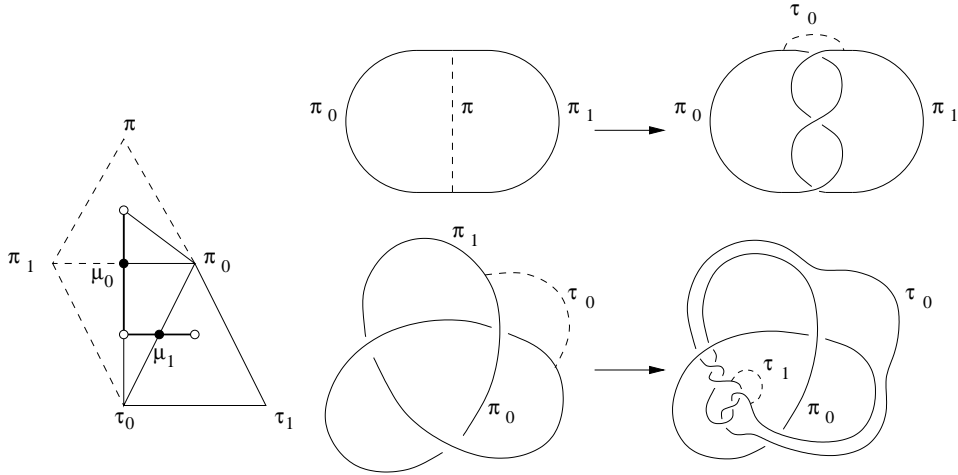}
\caption{Examples of the cabling construction.}
\label{fig:cabling}
\end{center}
\end{figure}

Figure~\ref{fig:cabling} shows the effects of a specific sequence of two
cabling constructions, starting with the trivial knot and obtaining the
trefoil, then starting with the tunnel of the trefoil.

As usual, let $\mu_0$ be a primitive pair, and let $\tau_0$ be a simple
disk for $\mu_0$, with simple slope $m_0$.  The segment $(\mu_0;\tau_0)$
determines a cabling construction starting with the tunnel of the trivial
knot and producing the $2$-bridge knot. We call this a \textit{simple
cabling of slope $m_0$.}

Now, consider a path of length $3$ in $\T$ determined by the four vertices
$\mu$, $\mu\cup \{\tau\}$, $\mu'$, $\mu'\cup\{\tau'\}$, that is, by a
succession from $(\mu;\tau)$ to $(\mu';\tau')$, as in
figure~\ref{fig:schematic}. Let $\sigma$ be the unique disk of $\mu-\mu'$
(the ``trailing'' disk), which in the case shown in
figure~\ref{fig:schematic} happens to be $\rho$.  Denoting by $m$ the
$(\mu';\sigma)$-slope of $\tau'$, we call the corresponding cabling
construction a \textit{cabling of slope~$m$.} We require that $\tau'\neq
\tau$, that is, cablings do not allow one to ``backtrack'' in $\T$. In
terms of slope, such a cabling would have~$m=\infty$.

On the other hand, any cabling construction in $H$ corresponds to an arc of
length $3$ starting at a black vertex in $\widetilde{\T}$, and any cabling
construction in $S^3$ that does not produce the trivial tunnel or a simple
tunnel corresponds to such an arc in $\T$. One of our main
theorems is now obvious; it is just a geometric restatement of the
Parameterization Theorem~\ref{thm:parameterization}.
\begin{uniquecablingsequencetheorem} 
Let $\tau$ be a tunnel of a nontrivial knot.
Let $\theta_0$, $\mu_0$, $\mu_0\cup \{\tau_0\}$, $\mu_1,\ldots\,$, $\mu_n$,
$\mu_n\cup \{\tau_n\}$ with $\tau_n=\tau$ be the principal path of $\tau$.
Then the sequence of $n+1$ cablings consisting of the simple cabling
determined by $(\mu_0;\tau_0)$ and the cablings determined by the
successions from $(\mu_{i-1};\tau_{i-1})$ to $(\mu_i;\tau_i)$ is the
\textit{unique} sequence of cablings beginning with the tunnel of the trivial
knot and ending with~$\tau$.
\label{thm:cable_sequence}
\end{uniquecablingsequencetheorem}
\noindent Of course, the values $[p_0/q_0]$ and $q_i/p_i$ in the
Parameterization Theorem~\ref{thm:parameterization} are exactly the slopes
of the cablings.

One of the slope parameters in the Unique Cabling Sequence
Theorem~\ref{thm:cable_sequence} will be the subject of
section~\ref{sec:STinvariant}:
\begin{definition} Let $\tau$ be a tunnel for a nontrivial knot.
If $\tau$ is not simple, then the slope $m_n$ is called the
\textit{principal slope} of~$\tau$. When $\tau$ is simple, its principal
slope is undefined.
\end{definition}

In some sense, theorem~\ref{thm:cable_sequence} enables one to distinguish
any two tunnels. If one finds any sequence of cablings that produce the
tunnel, it must be the \textit{unique} such sequence, and when the
sequences are different for two tunnels, the tunnels are inequivalent.  The
sequence of cablings, i.~e.~the principal path, can be determined
algorithmically starting from a specific representative disk $D\subset H$
of the tunnel: Working in $\widetilde{\T}$, fix a primitive pair
$\mu_0=\{\pi_0,\pi_1\}$, and let $W$ be a wave for $D$ with respect to
$\mu_0$ with $W$ meeting, say, $\pi_1$. Let $\tau_0$ be the slope disk of
$\mu_0$ determined by $W$. Then $D$ has fewer components of intersection
with $\pi_0\cup \tau_0$ than with $\pi_0\cup \pi_1$. Put
$\mu_1=\{\pi_0,\tau_0\}$ and repeat this process inductively using a wave
of $D$ with respect to $\mu_1$. When $D$ has no wave, one is at the
principal vertex $\mu_n\cup \{D\}$. Some of the initial disks $\tau_0$,
$\tau_1,\ldots\,$ may be primitive, but the portion starting with the last
$\mu_i\cup\{\tau_i\}$ that is a primitive triple will descend to the
principal path in $\T$.

\section{The Scharlemann-Thompson invariant}
\label{sec:STinvariant}

The Scharlemann-Thompson invariant, developed in
\cite{Scharlemann-Thompson} and further used in
\cite{ScharlemannGodaTeragaito, SaitoSTinvariant}, is essentially the
principal slope of $\tau$. The construction in \cite{Scharlemann-Thompson}
proceeds as follows: (1)~intersect a splitting sphere $S$ with $H$,
obtaining a separating disk $E$, (2)~take an arc of $E\cap \tau$ outermost
on $E$ and cutting off a subdisk $E'$ of $E$, (3)~take as $\mu^+$ and
$\mu^-$ the two components of the frontier of a regular neighborhood of
$\tau\cup E'$ that are not parallel to $\tau$.  Then the invariant is
defined to be the $(\{\mu^+,\mu^-\};\tau)$-slope of a wave of $S$. No
selection of a canonical perpendicular disk $\tau^0$ is made, so the
invariant is regarded as an element of $\Q/2\Z$ since changing the choice
of $\tau^\perp$ changes the slope by a multiple of $2$. We will regard the
Scharlemann-Thompson invariant as $\Q$-valued, by using $\tau^0$ as the
perpendicular disk.

\begin{figure}
\begin{center}
\includegraphics[width=48 ex]{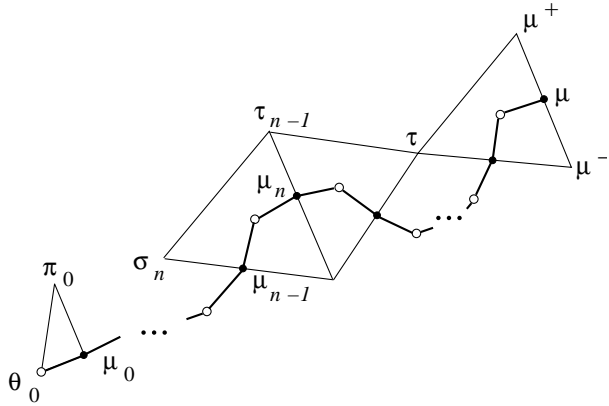}
\caption{Verification that the principal vertex is $\{\mu^+,\mu^-,\tau\}$.}
\label{fig:STcoords}
\end{center}
\end{figure}
To understand the relation between the Scharlemann-Thompson invariant and
the principal slope, we use the following lemma.
\begin{lemma} Let $\tau$ be a non-simple tunnel of a nontrivial knot. Then
\begin{enumerate}
\item[(i)] The disks called $\mu^+$ and $\mu^-$ in the definition of the
Scharlemann-Thompson invariant form the principal meridian pair of
$\tau$. That is, the principal vertex of $\tau$ is
$\{\tau,\mu^+,\mu^-\}$. 
\item[(ii)] The Scharlemann-Thompson invariant
is the $(\mu_n;\tau)$-slope of the unique disk $\sigma_n$ of
$\mu_{n-1}-\mu_n$.
\end{enumerate}
\label{lem:muplusminus}
\end{lemma}
\begin{proof}
Before beginning the proof, we remind the reader that any finite collection
of essential disks in $H$ can be moved by isotopies so that any two of them
intersect minimally. A nice way to do this is to choose a hyperbolic
structure on $\partial H$, move the boundaries of the disks by isotopies to
be geodesics, and eliminate simple closed curve intersections of their
interiors by futher isotopies. If we choose our hyperbolic structure so
that the hyperbolic involution $\alpha$ of $H$ is an isometry, then the
boundaries will be invariant under $\alpha$, indeed every simple closed
curve is invariant up to isotopy so every geodesic is invariant. In
particular, we may assume that $\alpha$ preserves the boundaries of $\tau$,
$\mu^+$, $\mu^-$, and~$E$.

Put $\mu=\{\mu^+,\mu^-\}$, and let $W$ be a wave of the separating disk $E$
with respect to $\mu$. We first claim that the $(\mu;\tau)$-slope of
$E$ is finite. If not, then $W$ is disjoint from $\tau$. Consider the arc
$W\cap \partial H$. Either this arc or its image under $\alpha$ lies on the
side of $\mu^+\cup \tau\cup \mu^-$ that contains $E'$; let $\omega$ be the
one that does. But $\omega$ is disjoint from $E'$, which contradicts the
fact that $W$ is essential.

Now, let $\mu_n\cup \{\tau\}$ be the principal vertex of $\tau$.
Figure~\ref{fig:STcoords} illustrates the path in $\T$ from $\theta_0$ to
$\mu$ if the principal vertex is not $\mu\cup\{\tau\}$ (it is possible that
one of $\mu^+$ or $\mu^-$ equals a disk of $\mu_n$).

The splitting sphere $S$ that contains $E$ is disjoint from a primitive
disk $D_0$. Suppose that $\mu\neq \mu_n$. Then we can travel through the
$1$-skeleton of $\D(H)$ from $D_0$ to $\tau$, passing through a sequence
$D_0$, $D_1,\ldots\,$, $D_k$ of vertices of $\D(H)$ such that:
\begin{enumerate}
\item For each $1\leq i\leq k$, $D_{i-1}$ is disjoint from $D_i$.
\item $D_{k-2}=\sigma_n$ and $D_k=\tau$.
\item No $D_i$ is a disk of $\mu$ (if it happens that one of $\mu^+$ or
$\mu^-$ is a disk of $\mu_n$, $D_{k-1}$ will be the other disk of $\mu_n$).
\end{enumerate}

Consider two disjoint, possibly separating disks in $H$, neither of which
is a disk of $\mu$. As before, they may be moved by isotopy to be disjoint
and simultaneously to have minimal intersection with $\mu$. If neither of
the two disks is a slope disk, then they have disjoint waves with respect
to $\mu$, which are disjoint from a unique slope disk of $\mu$. If one is a
slope disk, then it is disjoint from a wave of the other. In either case,
both are disjoint from the same slope disk of $\mu$, and since neither of
them was in $\mu$, this slope disk is uniquely determined. Inducting
through the sequence $E$, $D_0,\ldots\,$, $D_k=\tau$ now shows that a wave
of $E$ with respect to $\mu$ is disjoint from $\tau$, so the
$(\mu;\tau)$-slope of a wave of $E$ is infinite. This is a
contradiction.

For part~(ii), we have shown that $\mu=\mu_n$ and that a wave of any
splitting sphere is disjoint from~$D_{k-2}=\sigma_n$, so the
$(\mu_n;\tau)$-slope of a wave of $E$ is the same as the slope
of~$\sigma_n$.
\end{proof}

The previous argument clarifies the fact that the Scharle\-mann-Thompson
invariant is finite for exactly one meridian pair of $\tau$ (Lemma~2.9
of~\cite{Scharlemann-Thompson}). Also, it explains why the
Scharlemann-Thompson definition gives nonunique values for a few cases
(Corollary 2.8 of~\cite{Scharlemann-Thompson}). These cases are exactly
the upper and lower tunnels of $2$-bridge knots, for which the choice of
$\sigma_n$ is not unique. These are the simple tunnels, for which the
principal slope is undefined.

Theorem~\ref{thm:compare_invariants} below gives the expression for the
Scharlemann-Thompson invariant and the principal slope $m_n$ in terms of
each other. To obtain this, we must understand how to change coordinates on
the slope disks at the principal meridian pair of the tunnel, and we set
this up as a general principle. Recall that for integers $a_1,\ldots\,$,
$a_k$, the continued fraction $[a_1,\ldots,a_k]$ is defined inductively by
$[a_1]=a_1$ and $[a_1,\ldots,a_k] =a_1+1/[a_2,\ldots,a_k]$.  Sometimes one
may choose to allow some of the $a_i$ to be $\infty$. We have
$-[a_1,\ldots,a_k]=[-a_1,\ldots,-a_k]$, and
\[[\ldots a_{i-1},0,a_{i+1}\ldots]= 
[\ldots a_{i-1}+a_{i+1}\ldots]\]
so in particular
\[[\ldots,a_{k-1},a_k,0] =[\ldots,a_{k-1},a_k,0,\infty]
=[\ldots,a_{k-1},\infty] =[\ldots,a_{k-1}]\ .\]
Another basic fact is:
\begin{lemma} Let $q/p\in \Q$.
Then $q/p$ may be written as a continued fraction as
$[2a_1,2b_1,2a_2,\ldots,2b_{k-1},2a_k]$ or
$[2a_1,2b_1,2a_2,\ldots,2b_{k-1},2a_k,b_k]$, with all entries nonzero
except possibly $a_1$, according as $q$ is even or odd. When $q$ is odd,
the parity of $b_k$ equals the parity of~$p$. The expression is unique,
provided that $a_k$ and $b_k$ do not have different signs when $b_k=\pm1$.
\label{lem:even_continued_fraction}
\end{lemma}

There is a very well-known connection between continued fraction
decompositions and $\SL_2(\Z)$. Define
$U=\begin{pmatrix}1&1\\0&1\end{pmatrix}$ and
$L=\begin{pmatrix}1&0\\1&1\end{pmatrix}$.
\begin{lemma} If
$U^{a_1}L^{b_1}\cdots U^{a_k}L^{b_k} = \begin{pmatrix} q & s \\
p & r \end{pmatrix}$, then 
\begin{gather*}q/p = [a_1,b_1,\ldots, a_k, b_k],\;
s/r = [a_1,b_1,\ldots, a_k],\\
q/s = [b_k,a_k,\ldots, b_1, a_1]\text{, and }
p/r = [b_k,a_k,\ldots, b_1]\ .
\end{gather*}\par
\label{lem:continued_fractions}
\end{lemma}
\noindent The first two equalities can be proven by a straightforward
induction, and the last two follow by taking transposes.

\begin{proposition} Let $\mu=\{\lambda,\rho\}$ be a pair of disks in $H$ and
let $\sigma$ and $\tau$ be two nonseparating slope disks for $\mu$. Write
the $(\mu;\tau)$-slope of $\sigma$ as $q/p= [2a_1,2b_1,2a_2,
\ldots,2b_{n-1},2a_n,b_n]$. If we regard slope pairs $[p,q]$ as column
vectors $\begin{pmatrix} q \\ p\end{pmatrix}$, then the change-of-basis
matrix from $(\mu;\tau)$-slopes to $(\mu;\sigma)$-slopes is
\[U^{2a_1}L^{2b_1}U^{2a_2}\cdots U^{2a_n}L^{b_n}U^{-(-1)^p2a}\ ,\] 
where $a=\sum a_i$.
\label{prop:change_coordinates}
\end{proposition}
\begin{proof}
Referring to the picture of $H$ in figure~\ref{fig:slope_coords}, let $u$
be the Dehn twist of $H$ about $\tau$ that sends an object with slope pair
$[p,q]$ to one with slope pair $[p,q+2p]$ ($u$ is a ``left-handed'' Dehn
twist about $\tau$). Similarly, let $\ell$ be the homeomorphism of $H$ that
preserves $\lambda$ and $\rho$, and sends an object with slope pair $[p,q]$
to one with slope pair $[p+q,q]$ (a half-twist of $H$ about~$\tau^0$, whose
effect looks like the restriction of $\beta$ from section~\ref{sec:Akbas},
when $H$ is viewed as in figure~\ref{fig:slope_coords}).

Regarding slope pairs $[p,q]$ as column vectors $\begin{pmatrix} q \\
p\end{pmatrix}$, the effects of $u$ and $\ell$ are multiplication by $U^2$
and~$L$ respectively. That is, for the $(\mu;\tau)$-slopes determined by
the basis $\{\tau,\tau^0\}$, these are the matrices of $u$ and~$\ell$.

Since $\sigma$ is nonseparating, we can use
lemma~\ref{lem:even_continued_fraction} to write its
$(\mu;\tau)$-slope as $[2a_1,2b_1,\ldots,2b_{n-1},2a_n,b_n]$, where
$b_n$ has the parity of $p$ and all terms except possibly $2a_1$ are
nonzero. The composition $u^{a_1}\ell^{2b_1}\cdots u^{a_n}\ell^{b_n}$ of
$\ell$ has matrix $U^{2a_1}L^{2b_1}\cdots U^{2a_n}L^{b_n}$, so the
``$q/p$'' case of lemma~\ref{lem:continued_fractions} shows that it sends
$\tau$ to $\sigma$. It takes $\tau^0$ to a perpendicular disk for $\sigma$,
but not necessarily $\sigma^0$. Since $u$ preserves $\tau$, the
homeomorphism $u^{a_1}\ell^{2b_1}\cdots u^{a_n}\ell^{b_n}u^{-(-1)^pa}$,
where $a=\sum a_i$, also takes $\tau$ to $\sigma$. In the remainder of the
proof, we will show that it takes $\tau^0$ to~$\sigma^0$ as well, and the
lemma follows.

Consider any separating slope disk for $\mu$, and let $c_1$ and $c_2$ be
the core circles of its complementary solid tori in $H$. Each $c_i$ has
intersection number $0$ with $\tau^0$, but applying $\ell$ reverses the
sides of $\rho$, so reverses the orientation of exactly one of the
$c_i$. It follows that $\lk(\ell(c_1),\ell(c_2))=-\lk(c_1,c_2)$.  On the
other hand, each $c_i$ has intersection number $\pm 1$ with $\tau$, since
$c_i$ has intersection number $\pm 1$ with one of $\lambda$ or $\rho$ and
$0$ with the other. Therefore $\lk(u(c_1),u(c_2))=\lk(c_1,c_2)\pm 1$, the
sign depending on conventions.

In particular, if $C_1$ and $C_2$ are the core circles for the
complementary components of $\tau^0$, we have
\[\lk(\ell^{b_n}u^{-(-1)^pa}(C_1),\ell^{b_n}u^{-(-1)^pa}(C_2))
=\lk(u^{-a}(C_1)),u^{-a}(C_2))\ .\]
\noindent The remaining $\ell^{2b_i}$ do not change linking numbers, and
$a=\sum a_i$, so $\lk(u^{a_1}\ell^{2b_1}\cdots
u^{a_n}\ell^{b_n}u^{-(-1)^pa}(C_1), u^{a_1}\ell^{2b_1}\cdots
u^{a_n}\ell^{b_n}u^{-(-1)^pa}(C_2))=0$.  These are core circles of the
complementary tori of $u^{a_1}\ell^{2b_1}\cdots
u^{a_n}\ell^{b_n}u^{-(-1)^pa}(\tau^0)$. Since this disk intersects
$u^{a_1}\ell^{2b_1}\cdots u^{a_n}\ell^{b_n}u^{-(-1)^pa}(\tau)=\sigma$ in a
single arc, it must be~$\sigma^0$.
\end{proof}

\begin{theorem} If
$q/p = [2a_1,2b_1, 2a_2, \ldots, 2b_{n-1},2a_n,b_n]$ is one of the
Schar\-le\-mann-Thompson invariant or the principal slope of a tunnel, then
the other one is $[(-1)^p2a,-b_n,-2a_n,\ldots, -2a_2, -2b_1]$, where
$a=\sum a_i$. In particular, if one is a $($necessarily odd\/$)$ integer,
then the other is the negative of that integer.
\label{thm:compare_invariants}
\end{theorem}
\begin{proof} From lemma~\ref{lem:muplusminus},
the invariants are related by the fact that for some pair of slope disks
$\sigma$ and $\tau$ for a meridian pair $\mu$, one invariant is the
$(\mu;\sigma)$-slope of $\tau$, and the other is the $(\mu;\tau)$-slope of
$\sigma$. Using proposition~\ref{prop:change_coordinates}, the
change-of-basis matrix from $(\mu;\sigma)$-slopes to $(\mu;\tau)$-slopes is
$U^{(-1)^p2a}L^{-b_n}U^{-2a_n}\cdots U^{-2a_2}L^{-2b_1}U^{-2a_1}$. The
first column of this matrix gives the $(\mu;\tau)$-slope of $\sigma$. By
the ``$q/p$'' case of lemma~\ref{lem:continued_fractions}, with $b_k=0$, it
is $[(-1)^p2a,-b_n,-2a_n,\ldots, -2a_2, -2b_1]$.
\end{proof}

\begin{corollary} If $q/p$ is either the
Scharlemann-Thompson invariant or the principal slope of a tunnel, then the
other one is of the form $q\,'/p$ where $qq\,'\equiv -1\pmod{p}$.
\label{coro:convert}
\end{corollary}
\begin{proof} If the change-of-basis matrix in the proof of
theorem~\ref{thm:compare_invariants} has the form
$\begin{pmatrix}q&s\\p&r\end{pmatrix}$, then its inverse is
$\begin{pmatrix}r&-s\\-p&q\end{pmatrix}$, so the invariants are $q/p$ and
$-r/p$ where $qr-ps=1$.
\end{proof}

We have implemented the formula of theorem~\ref{thm:compare_invariants}
to convert between the invariants computationally~\cite{slopes}. Some
sample calculations are:
\smallskip

\noindent \texttt{STinvariant$>$ convert 55}\\
\texttt{-55}
\smallskip

\noindent \texttt{STinvariant$>$ convert (59/35)}\\
\texttt{-299/35}
\smallskip

\noindent \texttt{STinvariant$>$ convert (-299/35)}\\
\texttt{59/35}
\smallskip

\noindent \texttt{STinvariant$>$ convertRange 100102 17255 17265}\\
\texttt{17255/100102, -2843767/100102}\\
\texttt{17257/100102, -6541753/100102}\\
\texttt{17259/100102, 345051565/100102}\\
\texttt{17261/100102, 5593835/100102}\\
\texttt{17263/100102, 1775313/100102}\\
\texttt{17265/100102, 158447/100102}
\smallskip

\noindent The last command produces the corresponding pairs of invariants
containing each $q/100102$ for odd $q$ with $17255\leq q \leq 17265$.

\section{Tunnels of $2$-bridge knots}
\label{sec:2bridge}

It is known from work of Kobayashi \cite{Kobayashi1,Kobayashi2},
Morimoto-Sakuma~\cite{Morimoto-Sakuma}, and Uchida~\cite{Uchida}
that a $2$-bridge knot has at most four equivalence classes of tunnels (not
six, for us, since we are considering tunnels only up to equivalence,
rather than up to isotopy). Two of these are the upper and lower simple
tunnels. In this section, we will locate the other tunnels in~$\D(H)$ and
compute their slope parameters.

\begin{figure}
\begin{center}
\includegraphics[width=55 ex]{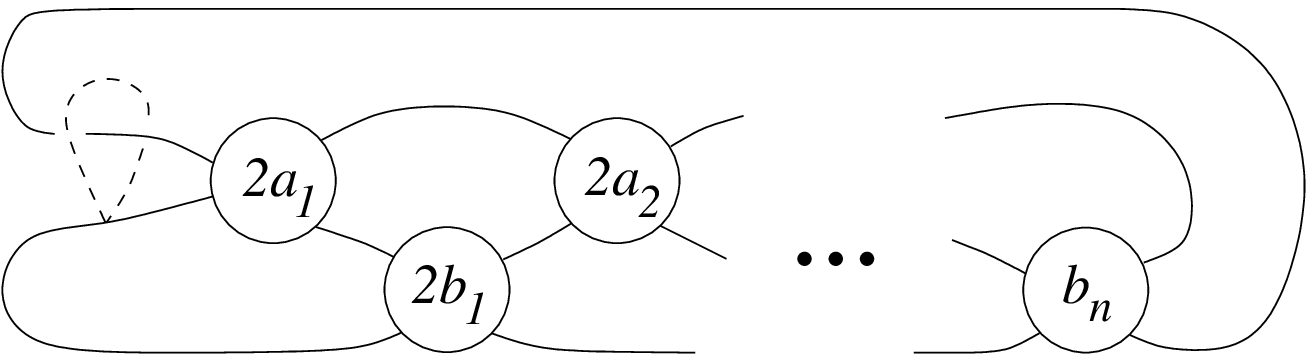}
\caption{}
\label{fig:2bridge}
\end{center}
\end{figure}
From the above references, and standard repositioning of $2$-bridge knots
by isotopy, each tunnel of a $2$-bridge knot either is simple or is
equivalent to one like that shown in figure~\ref{fig:2bridge}, where each
circle indicates a block of some nonzero number of half-twists. Each of the
blocks in the middle row has an even number $2a_i$ of half-twists, and
those on the bottom row have an even number $2b_i$ of half-twists, except
that the last one has a number $b_n$ that may be odd. Our convention is
that $a_i$ is positive for left-hand twists, and $b_i$ is positive for
right-hand twists.

There is a well-known classification of $2$-bridge knots based on continued
fraction expansions of a rational parameter $b/a$ with $b$ odd (the case of
$b$ even gives $2$-bridge links). One description of the invariant is that
the $2$-fold branched cover of $S^3$ over the knot is~$L(b,a)$, but we will
describe it here in a way that is more suited to our purposes.

Given $b/a$, change $a$ by multiples of $b$ until
$\vert\,b/a\,\vert>1$. Either of two possible values of $a$ may be
used. Expand $b/a$ as a continued fraction
$[2a_1,2b_1,2a_2,2b_2,\ldots,2a_n,b_n]$ as in
lemma~\ref{lem:even_continued_fraction}.  Additionally, if $b_n=\pm1$,
adjust $a_n$ and $b_n$ so that they have the same sign. Under these
conditions, the expansion of $b/a$ is uniquely determined, and the
corresponding $2$-bridge knot is the one shown in figure~\ref{fig:2bridge}.

\begin{figure}
\begin{center}
\includegraphics[width=\textwidth]{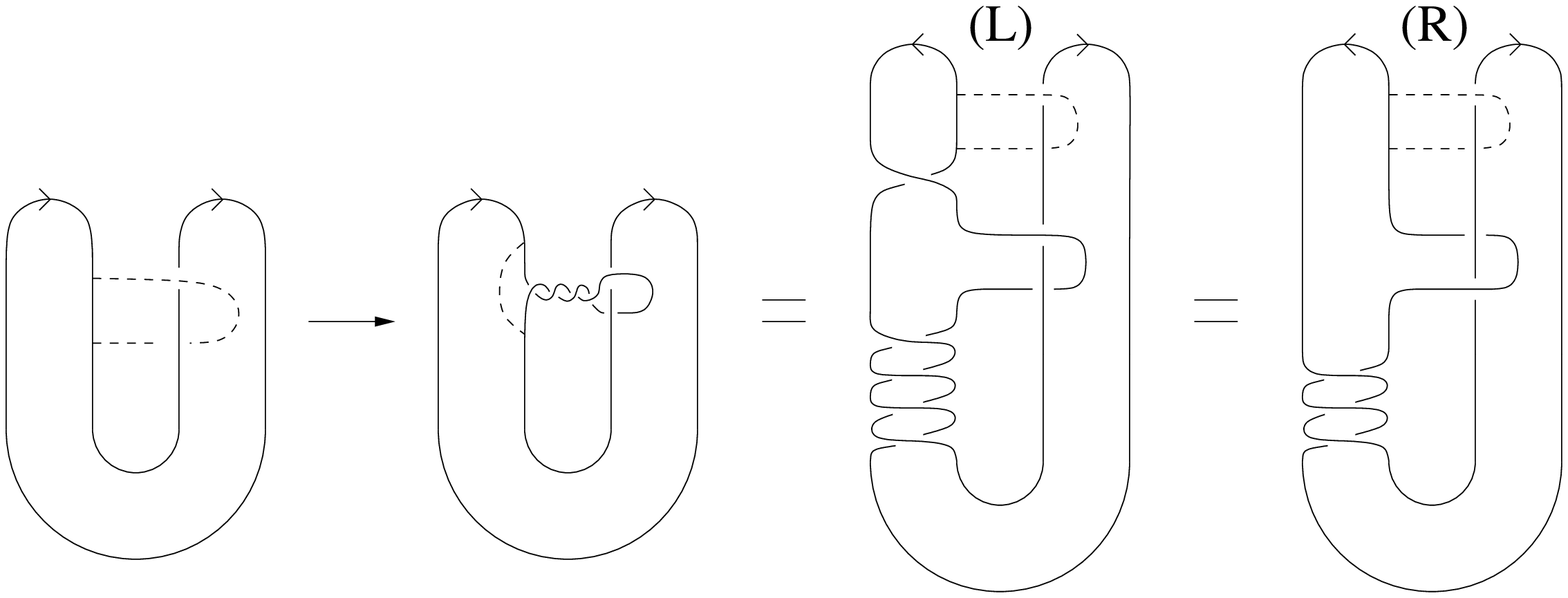}
\caption{}
\label{fig:2bridge_tunnel}
\end{center}
\end{figure}
Figure~\ref{fig:2bridge_tunnel} shows the type of cabling construction used
to produce the tunnel in figure~\ref{fig:2bridge}. It is described by a
nonzero parameter $k$ that tells the number of right-hand twists of the two
horizontal arcs from the original tunnel (the case $k=0$ would produce a
cabling with infinite slope, i.~e.~not a cabling construction). The cabling
shown in figure~\ref{fig:2bridge_tunnel} has $k=-4$. As indicated in that
figure, a knot and tunnel resulting from such a cabling can be moved by
isotopy so that the full twist of the middle two strands is either
left-handed (configuration~(L)) or right-handed (configuration~(R)), then
repositioned so that the tunnel has the same appearance as the original
one. As shown in figure~\ref{fig:2bridge_tunnel}, this will produce either
$k$ half-twists below the full twist and $-1$ half-twists above, as in
configuration~(L), or $k+1$ half-twists below the full twist, as in
configuration~(R).

Starting from the right-hand end of figure~\ref{fig:2bridge}, we perform a
sequence of these cablings, one for each of the full twists of the middle
two strands. Thus, the total number of cablings is $\sum\vert a_i\vert$. At
each step, the value of $k$ in the cablings must be selected to produce the
correct number $2b_i$ (or $b_n$) of half-twists, as we will detail
below. The condition that $b_n$ has the same sign as $a_n$, when
$b_n=\pm1$, ensures that the first cabling produces a nontrivial knot.  By
the Unique Cabling Sequence Theorem~\ref{thm:cable_sequence}, this sequence of
cablings is the unique sequence producing this tunnel.

We will now calculate the slopes of these cablings. The calculation depends
on the parity, at the time a given cabling is to be performed, of the
number of crossings of the left-hand two strands that lie below the tunnel
when it is positioned in figure~\ref{fig:2bridge_tunnel}.  In the example
of figure~\ref{fig:2bridge_tunnel}, the parity is odd both in
configuration~(L) and in configuration~(R). When the parity is even, either
of the two orientations of the knot orients the middle strands so that near
the tunnel one is upward and the other is downward, as occurs in the
trivial knot, but when the parity is odd, both are upward or both are
downward.

\begin{figure}
\begin{center}
\includegraphics[width=48 ex]{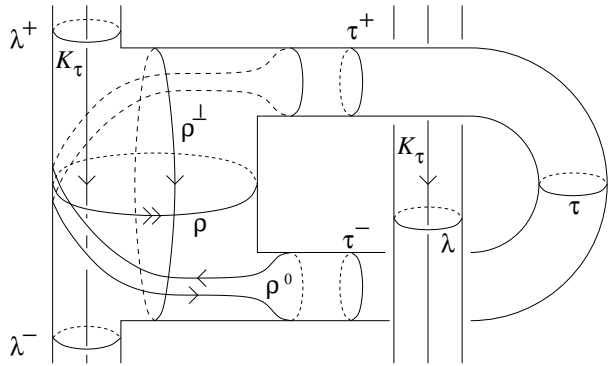}
\caption{The zero-slope disk for the case of odd parity
(compare with figure~\ref{fig:slope_coords}).}
\label{fig:odd_parity}
\end{center}
\end{figure}
The parity affects which disk will be the zero-slope disk in the
calculation of the slope coefficients. Figure~\ref{fig:odd_parity} shows
the case when the parity is odd. The disk $\rho$ is the one that is
replaced by the cabling construction, and the parity causes $\rho^0$ to be
as shown ($K_\lambda$ is a trivial knot that encircles the middle two
strands of $K_\tau$, and $\rho^0$ is obtained from $\rho^\perp$ by a Dehn
twist about $\rho$ that positions it as in the figure so as to make the
linking number of $K_\lambda$ and $K_\tau$ equal to $0$).

To obtain the actual slope coefficient for a given cabling of the type
shown in figure~\ref{fig:2bridge_tunnel}, one may draw a cabling arc for
the cabling and laboriously calculate its slope using $\rho^0$ as the
zero-slope disk, but there is a quick way to find it. If the disk
$\rho^\perp$ shown in figure~\ref{fig:odd_parity} were used in place of
$\rho^0$ in calculating the slope, then the slope pairs would be
$[k,1]$. The Dehn twist about $\rho$ that moves $\rho^0$ to $\rho^\perp$
moves the cabling arc to one whose slope pair using $\rho^\perp$ as the
zero-slope disk is $[k,1-2k]$, so this is the slope pair of the original
cabling arc using $\rho^0$. Consequently the slope of the cabling arc with
respect is $-2+1/k$. In the even parity case, the only difference is that
the twist moving $\rho^0$ to $\rho^\perp$ is in the opposite sense,
changing $[k,1]$ to $[k,1+2k]$ and producing a slope of~$2+1/k$.

Consider the initial cabling construction, which is an even-parity case.
Let $k$ be the number of right-hand half twists in the cabling, so that the
slope is $2+1/k$.

Suppose first that $a_n>0$. According to configuration~(L) of
figure~\ref{fig:2bridge_tunnel}, $b_n=k$. The parity of the trivial knot
is even, to the slope pair of the cabling is $[k,2k+1]=[b_n,2b_n+1]$,
giving $m_0=[b_n/(2b_n+1)]$.

Suppose now that $a_n<0$. From configuration~(R) of
figure~\ref{fig:2bridge_tunnel}, we have $b_n=k+1$ so the result is a
$2$-bridge knot with invariant $-2+1/(k+1)$. The slope pair of the cabling
is $[k,1+2k]=[b_n-1,2b_n-1]$, giving $m_0=[(b_n-1)/(2b_n-1)]$.

We remark that in either case, $m_0$ is represented by the reciprocal of
the standard invariant as we have described it here (in the latter case,
the reciprocal is $b_n/(-2b_n+1)$, and
$[-b_n/(2b_n-1)]=[(b_n-1)/(2b_n-1)]$ in $\Q/\Z$).

To examine the cablings beyond the first, it is notationally convenient to
rewrite the continued fraction by expanding each $2a_i$ to
$[2,0,2,0,\ldots,2,0,2]$ or $[-2,0,-2,0,\ldots,-2,0,-2]$, thereby assuming
that each $a_i=\pm2$ and allowing some $b_i=0$.

Consider the cabling that produces the full twist of the middle two strands
corresponding to $a_i$. The cabling that produces the twist corresponding
to $a_{i+1}$ has just been completed. Suppose first that $a_{i+1}>0$. As in
configuration~(L) of figure~\ref{fig:2bridge_tunnel}, there is already a
left-hand half twist in the left two strands. Since all $2b_j$ are
even, the parity is the opposite of the parity of $b_n$.  Again
referring to figure~\ref{fig:2bridge_tunnel}, we see that to end up with
exactly $b_i$ half-twists of the left two strands, we need to use $k=b_i+1$
if $a_i>0$, and $k=b_i$ if~$a_i<0$.

Suppose now that $a_{i+1}<0$. From configuration~(R) of
figure~\ref{fig:2bridge_tunnel}, the parity is just equal to that
of~$b_n$. To achieve $2b_i$ half-twists after the cabling, we need $k=2b_i$
if $a_i>0$, and $k=2b_i-1$ if~$a_i<0$.

We now have a complete algorithm to determine the values $k_i$ in the
cablings, and the cabling slopes: Write $b/a$ as
$[2a_1,2b_1,\ldots,2a_n,b_n]$ where
\begin{enumerate}
\item Each $a_i=\pm1$, and some $b_i$ other than $b_n$ may be $0$.
\item If $b_n=\pm1$, then $a_n$ and $b_n$ have the same sign.
\end{enumerate}
\noindent For the first cabling:
\begin{enumerate}
\item If $a_n=1$, then $k_n=b_n$ and $m_0=[b_n/(2b_n+1)]$.
\item If $a_n=-1$, then $k_n=b_n-1$ and $m_0=[(b_n-1)/(2b_n-1)]$.
\end{enumerate}
\noindent For the remaining cablings, the slope $m_i$ is $2+1/k_i$ or
$-2+1/k_i$ according as the parity is even or odd. The parity and the value
of $k_i$ are computed as follows:
\begin{enumerate}
\item If $a_{i+1}=1$, then:
\begin{enumerate}
\item The parity equals the parity of $b_n+1$.
\item If $a_i=1$, then $k_i=2b_i+1$, and if $a_i=-1$, then $k_i=2b_i$.
\end{enumerate}
\item If $a_{i+1}=-1$, then
\begin{enumerate}
\item The parity equals the parity of $b_n$.
\item If $a_i=1$, then $k_i=2b_i$, and if $a_i=-1$, then $k_i=2b_i-1$.
\end{enumerate}
\end{enumerate}

We have implemented the algorithm computationally~\cite{slopes}. Some
sample calculations are:
\smallskip

\noindent \texttt{TwoBridge$>$ slopes (33/19)}\\
\texttt{[ 1/3 ], 3, 5/3}
\smallskip

\noindent \texttt{TwoBridge$>$ slopes (64793/31710)}\\
\texttt{[ 2/3 ], -3/2, 3, 3, 3, 3, 3, 7/3, 3, 3, 3, 3, 49/24}
\smallskip

\noindent \texttt{TwoBridge$>$ slopes (3860981/2689048)}\\
\texttt{[ 13/27 ], 3, 3, 3, 5/3, 3, 7/3, 15/8, -5/3, -1, -3}
\smallskip

\noindent \texttt{TwoBridge$>$ slopes (5272967/2616517)}\\
\texttt{[ 5/9 ], 11/5, 21/10, -23/11, -131/66}
\medskip

Of course, the slope parameters that we have calculated are the parameters
$m_i$ that appear in the Parameterization
Theorem~\ref{thm:parameterization}. The parameters $s_j$ are all $0$,
since these are semisimple tunnels (the disk called $\lambda$ in
figure~\ref{fig:odd_parity} is retained in every cabling construction).

\section{Tunnels of links}
\label{sec:links}

A quick summary of how the theory adapts to include tunnels of tunnel
number $1$ links is that one just adds the separating disks as possible
slope disks. The cabling sequence ends with the first separating slope
disk, and cannot be continued. The Parameterization
Theorem~\ref{thm:parameterization} holds as stated, except allowing $q_n$
to be even.

In a bit more detail, one way to allow links is to include separating disks
in the theory from the start, that is, to use the full disk complex $\K(H)$
rather than the nonseparating disk complex $\D(H)$. Very little additional
complication actually occurs. A separating disk $E$ in $H$ is disjoint from
only two other disks, both nonseparating, so is a vertex of only one
$2$-simplex $\langle E,\tau_1,\tau_2\rangle$ attached to $\D(H)$ along
$\langle \tau_1,\tau_2\rangle$. Each $1$-simplex $\langle
\tau_1,\tau_2\rangle$ of $\D(H)$ is a face of countably many such
$2$-simplices, one for each separating slope disk of
$\{\tau_1,\tau_2\}$. In slope coordinates, these disks correspond to the
$q/p$ with $q$ even.

The link of $E$ in the first barycentric subdivision $\K'(H)$ consists of
two $1$-simplices meeting in the principal vertex $\{ E,\tau_1,\tau_2\}$
of~$E$. The spine of $\K(H)$ is obtained from $\widetilde{\T}$ simply by
adding a ``Y'' in each $2$-simplex $\langle E,\tau_1,\tau_2\rangle$, which
meets $\widetilde{\T}$ only in the vertex $\{\tau_1,\tau_2\}$.

A \textit{primitive} separating disk is a disk in $H$ that is contained in
a splitting sphere of $H$. Note that both of the nonseparating disks
disjoint from a primitive separating disk are primitive.

To obtain $\K(H)/\G$ from $\D(H)/\G$, we first add one half-simplex to the
primitive region. It meets the primitive simplex $\Pi$ along the edge
called $\langle \mu_0, \pi_0\rangle$ in section~\ref{sec:quotients}, and
its third vertex is the unique orbit of primitive separating disks. Next, a
half-simplex for each simple separating disk is added along $\langle \mu_0,
\pi_0\rangle$. These correspond to the $[p/q]\in \Q/\Z$ with $q$ even, and
are the upper and lower tunnels of $2$-bridge links. The remaining added
$2$-simplices are attached along the other $1$-simplices of $\D(H)/\G$ as
they were in~$\D(H)$. A tunnel of a $2$-bridge link has a principal path
from the primitive nonseparating triple $\theta_0$ to its principal vertex,
which is the only white vertex in the link of the tunnel.

The trivial link is the link associated to the orbit of primitive
separating disks. It arises from the tunnel of the trivial knot by a
cabling construction of simple slope~$[1/0] = \infty$. This is the only
case in which $\infty$ is an allowable slope parameter. With this convention,
we can state the general Parameterization Theorem:
\begin{theorem}
Let $\tau$ be a knot or link tunnel with principal path $\theta_0$, $\mu_0$,
$\mu_0\cup \{\tau_0\}$, $\mu_1,\ldots\,$, $\mu_n$, $\mu_n\cup \{\tau_n\}$.
Fix a lift of the principal path to $\K(H)$, so that each $\mu_i$
corresponds to an actual pair of disks in~$H$.
\begin{enumerate}
\item If $\tau$ is primitive, put $m_0=[0]\in\Q/\Z$ or $m_0=[1/0]=\infty$,
according as $\tau$ is the tunnel of the trivial knot or the trivial
link. Otherwise, let $m_0=[p_0/q_0]\in\Q/\Z$ be the simple slope of
$\tau_0$.
\item If $n\geq 1$, then for $1\leq i\leq n$ let $\sigma_i$ be the unique disk
in $\mu_{i-1}-\mu_i$ and let $m_i=q_i/p_i\in\Q$ be the
$(\mu_i;\sigma_i)$-slope of $\tau_i$.
\item If $n\geq 2$, then for $2\leq i\leq n$ define $s_i=0$ or $s_i=1$
according to whether or not the unique disk of $\mu_{i-1}\cap\mu_i$ equals the
unique disk of $\mu_{i-1}\cap\mu_{i-2}$.
\end{enumerate}
Then, sending $\tau$ to the pair $((m_0,\ldots,m_n),(s_2,\ldots,s_n))$ is a
bijection from the set of all tunnels of all tunnel number~$1$ knots to the
set of all elements $(([p_0/q_0],q_1/p_1,\ldots,q_n/p_n),(s_2,\ldots,s_n))$
in
\[\big(\Q/\Z\,\cup\,\{\infty\}\big) \,\cup\, 
\big(\Q/\Z\,\times\, \Q\big) \,\cup\, \big(\cup_{n\geq 2} \;\Q/\Z\,\times\,
\Q^n \,\times\,\, C_2^{n-1}\big)\] with all $q_i$ odd except possibly~$q_n$.
The tunnel is a tunnel of a knot or a link according as $q_n$ is odd or even.
\label{thm:link_parameterization}
\end{theorem}

The linking number of the two components of a tunnel number $1$ link, up to
sign, is half the numerator $q_n$ of the principal slope $q_n/p_n$ of
$\sigma$ (or half the denominator of the simple slope, if the tunnel is
simple). This is immediate from the construction of general slope coordinates in
section~\ref{sec:general_slopes}.

Theorem \ref{thm:link_parameterization} implies that a tunnel is almost
never equivalent to itself by an orientation-reversing equivalence:
\begin{theorem} Let $\tau$ be a tunnel of a tunnel number~$1$ knot or link.  
Suppose that $\tau$ is equivalent to itself by an orientation-reversing
equivalence. Then $\tau$ is the tunnel of the trivial knot, the trivial
link, or the Hopf link.
\label{thm:HopfLink}
\end{theorem}
\begin{proof}
As noted in remark~\ref{rem:ORequivalence}, the classification of tunnels
up to arbitrary homeomorphism of $S^3$ is obtained from that of
theorem~\ref{thm:parameterization} by adding the equivalence
$([m_0],m_1,\ldots, m_n)\sim ([-m_0],-m_1,\ldots, -m_n)$. The only tuples
equal to themselves under this move are $([0])$, $(\infty)$, and $([1/2])$,
which correspond to the trivial knot, the trivial link, and the Hopf link.
\end{proof}
\noindent The tunnel of the Hopf link will be examined more closely in
section~\ref{sec:Hopf_link} below.

The Unique Cabling Sequence Theorem~\ref{thm:cable_sequence} holds as
stated for links as well as knots.

As in proposition~\ref{prop:top_and_bottom}, the simple tunnels of links
are exactly the upper and lower tunnels of $2$-bridge links, and the
statement of proposition~\ref{prop:simple_classification} holds allowing
$q$ even. Concerning these tunnels, we can give a quick proof of a theorem
of C. Adams and A. Reid~\cite{Adams-Reid} and M. Kuhn~\cite{Kuhn}:
\begin{theorem}[Adams-Reid, Kuhn]
The only tunnels of a $2$-bridge link are its upper and lower tunnels.
\label{thm:2bridgelinks}
\end{theorem}
\begin{proof}
Since each component of a $2$-bridge link is unknotted, the tunnel disk
is disjoint from a primitive pair, hence is simple.
\end{proof}

We can also understand semisimple tunnels of links, that is, tunnels whose
cocore disk is disjoint from a primitive disk, but not from a primitive
pair. In this case, one of the components of the link is unknotted. Such
links are the topic of the following theorem, which slightly strengthens a
result of T. Harikae \cite{Harikae}:
\begin{theorem}
Let $L$ be a nontrivial tunnel number~$1$ link with an unknotted
component. Then the other component of $L$ is a $(1,1)$-knot. Moreover,
every tunnel of $L$ is simple or semisimple, and $L$ has torus bridge 
number~$2$.\par
\label{thm:semisimple_links}
\end{theorem}
\begin{proof}
Let $\sigma$ be a tunnel of $L$. Then $\sigma$ is disjoint from exactly two
nonseparating disks, $\tau_1$ and $\tau_2$, and one of them, say $\tau_1$,
must be primitive. Therefore $\tau_2$ is simple or semisimple, and
$K_{\tau_2}$ is a $(1,1)$-knot.

To prove that $L$ has torus bridge number~$2$, we refer to
figure~\ref{fig:semisimple_links}.
\begin{figure}
\begin{center}
\includegraphics[width=\textwidth]{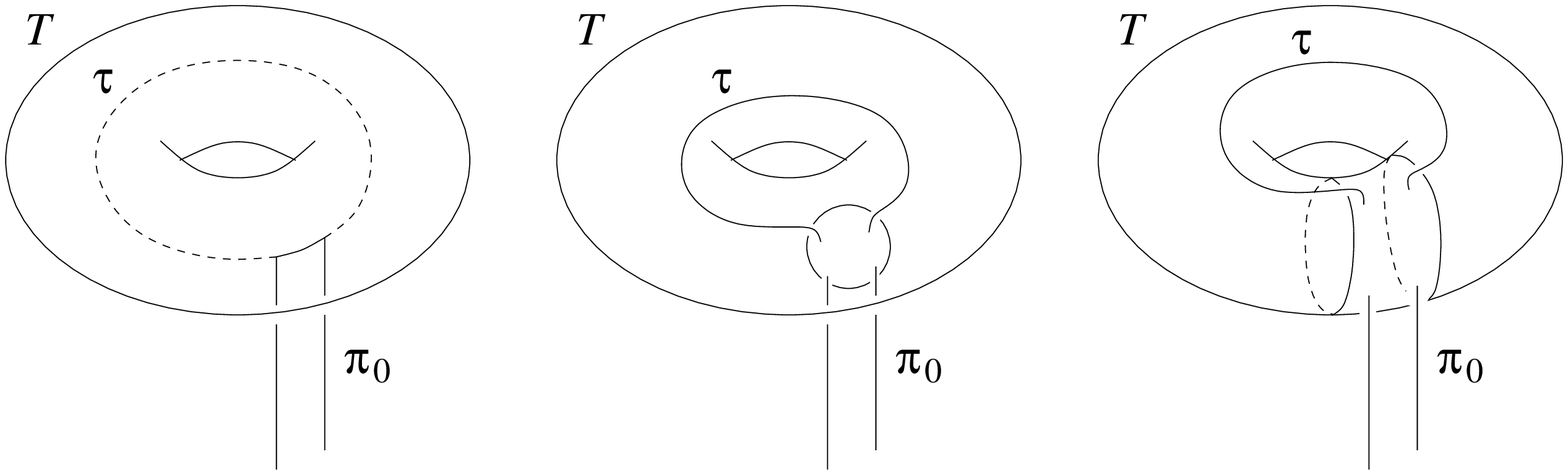}
\caption{}
\label{fig:semisimple_links}
\end{center}
\end{figure}
The left drawing shows a torus level $T$ and the tunnel $\tau$ before the
cabling that produces $\sigma$, and the middle picture shows schematically
the result of the cabling.  The drawing on the right shows an isotopic
repositioning of $L_\sigma$. The $\tau$-arc is pushed slightly outside of
$T$, and the sphere for the cabling is expanded to the union of an annulus
in $T$ and two meridian disks. The cabling arcs may be moved off of the
meridian disks by isotopy, to lie in $A$. From there, they can be pushed
slightly inside $T$. Then, $L$ meets each of the complementary solid tori
of $T$ in a trivial pair of arcs.
\end{proof}

\section{The Hopf link}
\label{sec:Hopf_link}

In theorem~\ref{thm:HopfLink}, we saw that the tunnel of the Hopf link is
the only tunnel of a nontrivial knot or link that is preserved by an
orientation-reversing equivalence. In this section, we examine this
equivalence more closely, obtaining a version of Futer's
Theorem~\ref{thm:Futer} for links.

First, we work out the vertex stabilizers of the action of $\G_\pm$ on
$\K(H)$. As usual, let $\mu_0$ be the standard primitive pair.  To obtain a
generating set for the stabilizer $(\G_\pm)_{\mu_0}$, we add to the
generators $\alpha$, $\beta$, and $\gamma$ of $\G_{\mu_0}$ an
orientation-reversing involution $R$ that commmutes with $\alpha$ and
$\gamma$ and conjugates $\beta$ to $\beta^{-1}$ (in the standard picture of
$H$, $R$ reflects in the plane of the page). Observe that $R$ sends a
simple disk of slope pair $[p,q]$ to one with slope pair $[-p,q]$. We have
the following version of proposition~\ref{prop:stabilizer}:
\begin{proposition} The stabilizer $(\mathcal{G}_\pm)_{\mu_0}$ is the subgroup
generated by $\alpha$, $\beta$, $\gamma$, and $R$. In fact, $(\G_\pm)_{\mu_0}$
is the semidirect product $(C_2\times \Z)\circ (C_2\times C_2)$, where
$\langle \alpha,\beta\rangle$ is the normal subgroup $C_2\times \Z$,
$\langle \gamma, R\rangle$ is the subgroup $C_2\times C_2$, $\alpha$ is
central, $\gamma\beta\gamma^{-1}=\alpha\beta$, and $R\beta
R^{-1}=\beta^{-1}$.
\label{prop:extended_stabilizer}
\end{proposition}

\begin{proposition} The stabilizer in $\G_\pm$
of a possibly separating nonprimitive disk $E$ in $H$ is as follows:
\begin{enumerate}
\item[(i)] If $E$ is simple with simple slope $[1,2]$, then its stabilizer is
conjugate to the dihedral subgroup of order $8$ generated by the involutions
$\beta R$ and~$\gamma$.
\item[(ii)] If $E$ is simple with simple slope not $[1,2]$, then its stabilizer
is conjugate to the subgroup $C_2\times C_2$ generated by $\alpha$ and~$\gamma$.
\item[(iii)] Otherwise, the stabilizer is $C_2$, generated by $\alpha$.
\end{enumerate}
\label{prop:extended_tau_stabilizer}
\end{proposition}

\begin{figure}
\begin{center}
\includegraphics[width=0.9\textwidth]{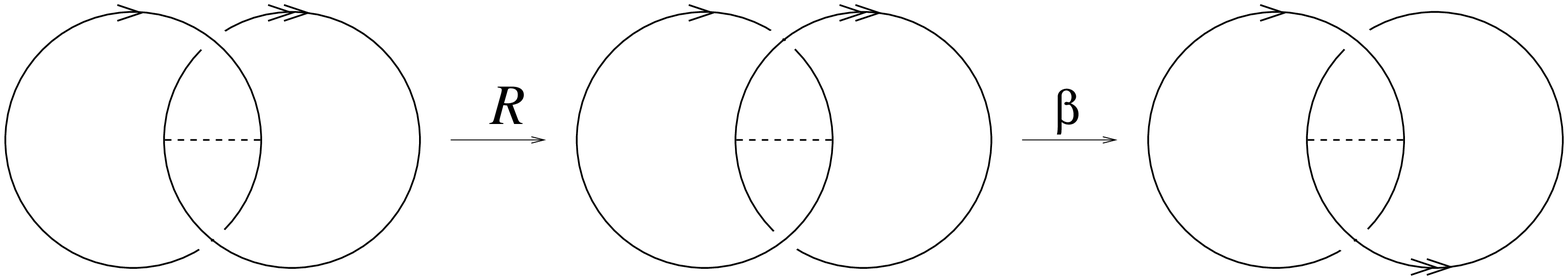}
\caption{An orientation-reversing automorphism of the Hopf link and its
tunnel arc.}
\label{fig:Hopf}
\end{center}
\end{figure}

\begin{proof} Assume first that $E$ is simple with simple slope $[1/2]$. 
Any element of $\G_\pm$ that preserves $E$ must also preserve the unique
pair of primitive disks that are disjoint from $E$. Conjugating in
$\G_\pm$, we may assume that this pair is the standard primitive pair
$\mu_0$.

Conjugating further by a power of $\beta$, we may assume that $E$ has slope
pair $[1,2]$ with respect to $\mu_0$. This slope pair is preserved by
$\beta R$, $\gamma$, and $\alpha$. On the other hand, using the action of
the four generators of $(\G_\pm)_{\mu_0}$ on slope pairs, any word that
stabilizes $E$ can be written as
$\alpha^{\epsilon_1}\gamma^{\epsilon_2}(\beta R)^{\epsilon_3}$ where the
$\epsilon_i\in\{0,1\}$. Since $(\gamma \beta R)^2=\alpha$, the stabilizer
is as in~(i).

If $E$ is simple but does not have simple slope $[1,2]$, then $R$ does not
preserve its simple slope. So any element $\alpha^{\epsilon_1}\beta^n
\gamma^{\epsilon_2} R^{\epsilon_3}$ stabilizing $E$ must have
$\epsilon_3=0$. As in proposition~\ref{prop:tau_stabilizer}(i), part~(ii)
follows. Part~(iii) is proven as in
proposition~\ref{prop:tau_stabilizer}(ii); note that the disk called
$\tau_0$ there must be nonseparating, since otherwise the path could not
continue on to~$\theta$.
\end{proof}
\noindent Figure~\ref{fig:Hopf} illustrates the effect of $\beta R$.

\newpage
Arguing as in corollary~\ref{coro:simplex_stabilizers}, we have
\begin{corollary} Let $\langle \tau_1,\tau_2,\tau_3\rangle$ be a
$2$-simplex of $\D(H)$, with $\tau_3$ nonprimitive. Suppose that $\langle
\tau_1,\tau_2,\tau_3\rangle$ is stabilized by an element $h$ of $\G_\pm$
other than the identity or $\alpha$. Then $\tau_3$ is simple, and $\tau_1$
and $\tau_2$ are primitive.
\label{coro:extended_simplex_stabilizers}
\end{corollary}

Using corollary~\ref{coro:extended_simplex_stabilizers}, and
proposition~\ref{prop:extended_tau_stabilizer}, we adapt
corollary~\ref{coro:theta_homeos} to links:
\begin{corollary} Let $W$ be the union of a nontrivial tunnel 
number~1 link $L$ and a tunnel arc for $L$, and suppose that $h$ is a
homeomorphism of $S^3$ that preserves $W$. Then either
\begin{enumerate}
\item[(i)] $h$ is isotopic preserving $W$ to a homeomorphism which is the
identity or the hyperelliptic involution on a neighborhood of $W$, or
\item[(ii)] $L$ is a $2$-bridge link.
\end{enumerate}
Moreover, if $h$ is orientation-reversing, then $L$ is the Hopf link.
\label{coro:W_homeos}
\end{corollary}

Using these results, we obtain a version of Futer's Theorem~\ref{thm:Futer}
for links:
\begin{theorem} Let $A$ be a tunnel arc for a nontrivial link
$L\subset S^3$. Then:
\begin{enumerate}
\item[(a)] There exists an orientation-preserving homeomorphism of $S^3$ that
preserves $L\cup A$ and interchanges the components of $L$ if and only if
$L$ is a two-bridge link.
\item[(b)] There exists an orientation-reversing homeomorphim of $S^3$ that
preserves $L\cup A$ if and only if $L$ is the Hopf link.
\end{enumerate}
\label{thm:Futer_for_links}
\end{theorem}

\begin{proof} For a $2$-bridge link, $\gamma$ interchanges the components 
and preserves the tunnel arc, up to isotopy, so
corollary~\ref{coro:W_homeos} gives~(a). Part~(b) follows easily from
proposition~\ref{prop:extended_tau_stabilizer}.
\end{proof}

\bibliographystyle{amsplain}

\begin{thebibliography}{10}

\bibitem{Adams-Reid} C. Adams, A. Reid,
Unknotting tunnels in two-bridge knot and link complements,
\textit{Comment. Math. Helv.} 71 (1996), 617--627.

\bibitem{Akbas} E. Akbas, A presentation of the automorphisms of the
$3$-sphere that preserve a genus two Heegaard splitting, Mathematics ArXiv
math.GT/0504519.

\bibitem{B-R-Z} M. Boileau, M. Rost, and H. Zieschang, On Heegaard
decompositions of torus knot exteriors and related Seifert fibre spaces,
\textit{Math. Ann.} 279 (1988), 553--581.

\bibitem{Cho} S. Cho, Homeomorphisms of the $3$-sphere that preserve a
genus~$2$ Heegaard splitting, Mathematics ArXiv math.GT/0611767, to appear
in \textit{Proc. Amer. Math. Soc.}

\bibitem{slopes} S.Cho and D. McCullough, software
available at 
\texttt{www.math.ou.edu/$_{\widetilde{\phantom{n}}}$dmccullough/}~.

\bibitem{CM} S.Cho and D. McCullough, The depth invariant of knot tunnels,
preprint. 

\bibitem{Futer} D. Futer, Involutions of knots that fix unknotting tunnels,
\textit{J. Knot Theory Ramifications} 16 (2007), 741-748.

\bibitem{G-S-T} H. Goda, M. Scharlemann, A. Thompson, Levelling an
unknotting tunnel, \textit{Geom. Topol.} 4 (2000), 243--275.

\bibitem{Goeritz} L. Goeritz, Die Abbildungen der Brezelfl\"ache und der
Volbrezel vom Gesschlect $2$, \textit{Abh. Math. Sem. Univ. Hamburg} 9
(1933) 244--259.

\bibitem{Harikae} T. Harikae, On the triviality of bouquets and tunnel
number one links, \textit{Proceedings of the Workshop on Graph Theory and
Related Topics (Sendai, 1999),} \textit{Interdiscip. Inform. Sci.} 7
(2001), 1--3.

\bibitem{JohnsonBridgeNumber} J. Johnson, Bridge number and the curve
complex, Mathe\-mat\-ics Ar\-Xiv math.GT/\allowbreak0603102.

\bibitem{Johnson-Thompson} J. Johnson and A. Thompson, On tunnel number one
knots which are not $(1,n)$, Mathematics ArXiv math.GT/0606226.

\bibitem{JDiff} D. McCullough, Virtually geometrically finite mapping class
groups of 3-manifolds, \textit{J. Diff. Geom.} 33 (1991), 1--65.

\bibitem{Kobayashi1} T. Kobayashi, A criterion for detecting inequivalent
tunnels for a knot, \textit{Math. Proc. Cambridge Philos. Soc.} 107 (1990),
483--491.

\bibitem{Kobayashi2} T. Kobayashi, Classification of unknotting tunnels for
two bridge knots, in \textit{Proceedings of the Kirbyfest} (Berkeley, CA,
1998), 259--290 (electronic), Geom. Topol. Monogr., 2, Geom. Topol. Publ.,
Coventry, 1999.

\bibitem{Kuhn} M. Kuhn, Tunnels of $2$-bridge links, \textit{J. Knot Theory
Ramifications} 5 (1996), 167--171.

\bibitem{MMS} Y. Minsky, Y. Moriah, and S. Schleimer, High distance knots,
Mathematics ArXiv math.GT/0607265.

\bibitem{Morimoto-Sakuma} K. Morimoto, M. Sakuma, On unknotting
tunnels for knots, \textit{Math. Ann.} 289 (1991), 143--167.

\bibitem{ScharlemannTree} M. Scharlemann, Automorphisms of the 3-sphere
that preserve a genus two Heegaard splitting,
\textit{Bol. Soc. Mat. Mexicana} (3) 10 (2004) 503--514.

\bibitem{ScharlemannGodaTeragaito} M. Scharlemann, There are no unexpected
tunnel number one knots of genus one, \textit{Trans. Amer. Math. Soc.}  356
(2004), 1385--1442.

\bibitem{Scharlemann-Thompson} M. Scharlemann and A. Thompson, Unknotting
tunnels and Seifert surfaces, \textit{Proc. London Math. Soc.} (3) 87
(2003), 523--544.

\bibitem{Scharlemann-Tomova} M. Scharlemann and M. Tomova, Alternate
Heegaard genus bounds distance, \textit{Geom. Topol.} 10 (2006), 593--617.

\bibitem{SaitoSTinvariant} T. Saito, Scharlemann-Thompson
invariant for knots with unknotting tunnels and the distance of
$(1,1)$-splittings, \textit{J. London Math. Soc.} (2) 71 (2005), 801--816.

\bibitem{Uchida} Y. Uchida, Detecting inequivalence of some unknotting
tunnels for two-bridge knots, \textit{Algebra and Topology 1990,} Korea
Adv. Inst. Sci. Tech. Taejon (1990) 227--232.

\end{thebibliography}

\end{document}